\documentclass[11pt,twoside]{article}
\usepackage[hmargin=0.8in,vmargin=0.9in]{geometry}
\usepackage{fancyhdr}
\usepackage{graphicx}
\usepackage{subfigure}
\usepackage{amssymb}
\usepackage{amsmath}
\usepackage{amsfonts}
\usepackage{theorem}
\usepackage{mathrsfs}
\usepackage{mathtools}
\usepackage{bm}
\usepackage{color}
\usepackage{setspace}
\usepackage{exscale}
\usepackage{hyperref}
\usepackage{relsize}
\usepackage{epstopdf}
\usepackage{float}
\usepackage{cite}

\usepackage{algorithm}
\usepackage{algpseudocode}

\usepackage{color}

\usepackage{hyperref}
\usepackage{setspace}
\usepackage{placeins}

\usepackage{booktabs,multirow} 
\usepackage{array} 
\usepackage{paralist} 
\usepackage{verbatim} 
\usepackage{subfigure} 

\theoremstyle{plain}			

{\theorembodyfont{\rmfamily}\newtheorem{remark}{Remark}[section]}

\newenvironment{DA}{{\flushleft \bf Declarations:}}{}

\allowdisplaybreaks[1]

\numberwithin{equation}{section}
\numberwithin{figure}{section}
\numberwithin{table}{section}

\newcommand\eref[1]{(\ref{#1})}

\newcommand*\xbar[1]{%
  \hbox{%
    \vbox{%
      \hrule height 0.5pt 
      \kern0.4ex
      \hbox{%
        \kern-0.05em
        \ensuremath{#1}%
        \kern-0.00em
      }%
    }%
  }%
}

\newcommand{\mF}{\bm{F}}

\newcommand{\mG}{\bm{G}}

\newcommand{\mU}{\bm{U}}

\newcommand{\mo}{\bm{0}}

\newcommand{\dt}{\Delta t}
\newcommand{\dx}{\Delta x}
\newcommand{\dy}{\Delta y}
\newcommand{\eps}{\varepsilon}

\newcommand{\hf}{{\frac{1}{2}}}

\newcommand{\jph}{{j+\frac{1}{2}}}
\newcommand{\jmh}{{j-\frac{1}{2}}}
\newcommand{\kph}{{k+\frac{1}{2}}}
\newcommand{\kmh}{{k-\frac{1}{2}}}

\title{Adaptive Artificial Anti-Diffusion Methods for Hyperbolic Systems of Conservation Laws}
\author{Shaoshuai Chu\thanks{Department of Mathematics, RWTH Aachen University, Aachen, 52056, Germany; {\tt chu@igpm.rwth-aachen.de}},
~Igor Kliakhandler\thanks{Department of Mathematics, Michigan Technological University, Houghton, MI 49931, USA; {\tt igor@mtu.edu}}, ~and
Alexander Kurganov\thanks{Department of Mathematics and Shenzhen International Center for Mathematics, Southern University of Science and
Technology, Shenzhen, 518055, China; {\tt alexander@sustech.edu.cn}}}
\date{}
\begin{document}
\date{}
\maketitle

\begin{abstract}
We introduce new adaptive artificial anti-diffusion (AAAD) methods for one- and two-dimensional hyperbolic systems of conservation laws. The
key idea is to reduce the amount of numerical dissipation present in a given numerical method by adding an anti-diffusion (AD) term acting
in linearly degenerate fields only. This way, the resolution of contact waves can be improved without risking oscillations, which may be
caused if the AD acts in nonlinear fields as well. The AD coefficients are selected adaptively: they are supposed to be proportional to the
mesh size near the contact waves to enhance the resolution and to be very small in the smooth parts of the computed solution to ensure a
sufficiently high (formal) order of accuracy there. The proposed AAAD methods are realized using either the second-order central-upwind
numerical fluxes or their fifth-order extensions implemented within the alternative weighted essentially non-oscillatory (A-WENO) framework.
We test the proposed schemes on a series of benchmarks for the one- and two-dimensional Euler equations of gas dynamics and the obtained
results demonstrate the robustness and high resolution of the new AAAD methods.
\end{abstract}

\noindent
{\bf Keywords:} Adaptive artificial anti-diffusion; local characteristic decomposition; linearly degenerate fields; Euler equations of gas
dynamics; central-upwind schemes; A-WENO schemes.

\smallskip
\noindent
{\bf AMS subject classification:} 65M08, 65M06, 76M12, 76M20, 76L05, 35L65.

\section{Introduction}
This paper focuses on numerical solutions of hyperbolic systems of conservation laws, which in the one-dimensional (1-D) and
two-dimensional (2-D) cases, read as
\begin{equation}
\mU_t+\mF(\mU)_x=\mo,
\label{1.1}
\end{equation}
and
\begin{equation}
\mU_t+\mF(\mU)_x+\mG(\mU)_y=\mo,
\label{1.2}
\end{equation}
where $t$ denotes time, $x$ and $y$ are spatial coordinates, $\mU\in\mathbb R^d$ is the vector of conserved variables, and
$\mF,\mG:\mathbb R^d\to\mathbb R^d$ are the flux functions.

It is well-known that even for smooth initial data, solutions of \eref{1.1} and \eref{1.2} may develop complicated wave patterns containing
shock waves, rarefaction waves, and contact discontinuities. The presence of such nonsmooth structures makes the design of accurate and
robust numerical methods particularly challenging. On the one hand, a good scheme has to suppress nonphysical oscillations near
discontinuities in order to remain stable. On the other hand, it must retain sufficient resolution so that physically relevant structures
are not over-smeared by excessive numerical dissipation. Achieving this balance between robustness and low numerical dissipation is one of
the central difficulties in the development of numerical methods for hyperbolic conservation laws.

A classical approach to enforce stability is to add an artificial viscosity (AV) to the underlying unstable (or, only linearly stable)
numerical scheme. Since its introduction in \cite{NR1950}, this idea has been successfully used in many works; see, e.g.,
\cite{CSW1998,GP2008,Kurganov12a,HH2002,HH2002a,Szepessy1989,Wilkins1980} to name just a few. The main advantage of AV methods lies in their
simplicity and robustness, since they provide an effective mechanism for suppressing spurious oscillations and enhancing nonlinear
stability. However, designing a successful AV method remains a delicate task. On the one hand, a sufficient amount of AV must be added near
discontinuities in order to suppress spurious oscillations and guarantee nonlinear stability. If it is too small, the added AV may be
insufficient to stabilize the computation near discontinuities, and nonphysical oscillations may still occur. On the other hand, the added
AV should remain sufficiently small in smooth regions so that the high-order accuracy of the underlying scheme is not destroyed. If it is
too large, the numerical solution may become excessively diffused, leading to a substantial loss of resolution. Therefore, in order to
achieve both robustness and high resolution, the AV must be introduced adaptively, based on a suitable smoothness indicator (SI) capable of
detecting ``rough'' parts of the computed solution and determining the appropriate amount of AV to be added there. This has been done, for
instance, in \cite{GP2008,Kurganov12a}, but even these advanced adaptive artificial viscosity (AAV) methods do not achieve the same level of
resolution achieved by modern shock-capturing methods. 

In  high-resolution shock-capturing methods, stability is achieved by incorporating nonlinear limiters into the reconstruction or flux
evaluation procedures. This strategy has proved to be highly effective and underlies many robust non-oscillatory methods. However, the
robustness gained in this way is often accompanied by increased numerical dissipation. In particular, contact discontinuities and
small-scale smooth structures may become noticeably smeared, especially on relatively coarse meshes. Therefore, for numerical methods that
already employ nonlinear limiters, it is natural to seek a mechanism that improves resolution without destroying a non-oscillatory property
ensured by the limiter.

In this paper, we develop a new class of adaptive artificial anti-diffusion (AAAD) methods for hyperbolic systems of conservation laws. Our
starting point is somewhat similar to the AAV philosophy, but rather than locally increasing the amount of numerical diffusion, we
adaptively add suitable anti-diffusive (AD) corrections to enhance the resolution of the scheme, which is already stabilized by nonlinear
limiters. In this way, the nonlinear limiters continue to control spurious oscillations, while the AD terms compensate excessive numerical
diffusion. This is particularly advantageous near contact discontinuities, whose accurate resolution is often degraded by the numerical
diffusion inherent in shock-capturing schemes.

A key point in our construction is that the AD correction is applied in linearly degenerate characteristic fields only. This allows us to
sharpen contact waves without introducing undesirable oscillations into the genuinely nonlinear fields associated with shocks and
rarefactions. Another essential feature of the proposed AAAD method is adaptivity. The AD coefficients are chosen so that the correction is
sufficiently strong only in ``rough'' regions, while remaining very small in smooth parts of the computed solution. As a result, the
designed (formal) order of accuracy of the underlying scheme is retained in smooth regions, whereas the resolution of contact
discontinuities is significantly enhanced.

To adjust the AAAD coefficients, we need to automatically detect the ``rough'' parts of the computed solution with the help of an SI. Many
different SIs are readily available; see, e.g.,
\cite{DZLD14,Dewar15,ABD08,GLHWDW,GT06,GT02,GPP,PupSem,QS05,VR16,FS17,WSYK15,KKP02,RL87,WDGH20,CCK23_Adaptive} and references therein. In
this paper, we employ the recently proposed modified minmod-based SI introduced in \cite{CHK2026}, which is a modified version of the
minmod-based SI from \cite{WSYK15}; see also \cite{Harten89,SO89}. We apply this SI to different fields in order to automatically
distinguish contact waves from other ``rough'' parts of the solution. For example, in the case of the Euler equations of gas
dynamics, we take advantage of the fact that, while the density is discontinuous across contact waves, the pressure remains continuous
there. After identifying the locations of the contact surfaces, we activate the AAAD correction in their vicinities with coefficients of
order ${\cal O}(\Delta)$, while taking much smaller adaption coefficients elsewhere: ${\cal O}(\Delta^2)$ and ${\cal O}(\Delta^r)$ in the
rest of the ``rough'' and ``smooth'' regions, respectively. Here, $\Delta$ is a spatial mesh size and $r$ is a formal order of the
underlying scheme so that the resulting AAAD method retains its formal order of accuracy in ``smooth'' regions.

We incorporate the proposed AAAD strategy into both second- and fifth-order schemes. We stress that the AAAD technique is not tied to any
particular underlying numerical methods. We demonstrate the performance and potential advantages of the proposed AAAD methods using the
central-upwind (CU) numerical fluxes. First, we combine the AAAD with the second-order finite-volume (FV) semi-discrete CU schemes developed
in \cite{Kurganov01,Kurganov00,Kurganov07,Kurganov02}. Second, we extend the same idea to the fifth-order finite-difference (FD) alternative
weighted essentially non-oscillatory (A-WENO) framework developed in \cite{JSZ,Liu17,Liu16,Wang18,CKX22}. The resulting adaptive schemes are
then applied to a number of 1-D and 2-D examples for the Euler equations of gas dynamics. The reported numerical results demonstrate that
the proposed AAAD methods are robust and substantially less dissipative than their non-adaptive counterparts. In particular, they provide a
sharper resolution of contact waves and finer structures while preserving the non-oscillatory character of the underlying schemes.

The rest of the paper is organized as follows. In \S\ref{sec2}, we introduce the 1-D second- and fifth-order AAAD methods, and then extend
them to the 2-D case in \S\ref{sec3}. In \S\ref{sec4}, we present a number of 1-D and 2-D numerical examples illustrating the performance of
the proposed AAAD methods. Finally, in \S\ref{sec5}, we provide concluding remarks.

\section{One-Dimensional Adaptive Artificial Anti-Diffusion Schemes}\label{sec2}
In this section, we present the proposed second- and fifth-order AAAD schemes for the 1-D Euler equations of gas dynamics, which read as
\eref{1.1} with
\begin{equation}
\bm U=(\rho,\rho u,E)^\top,\quad\bm F(\bm U)=(\rho u,\rho u^2+p,u(E+p))^\top,
\label{2.3}
\end{equation}
where $\rho$, $u$, $p$, and $E$ denote the density, velocity, pressure, and total energy, respectively. The system is closed by the equation
of state (EOS) for ideal gases:
\begin{equation}
p=(\gamma-1)\Big[E-\hf\rho u^2\Big],
\label{2.4}
\end{equation}
where $\gamma$ is the ratio of specific heats.

The Jacobian of the system \eref{1.1}, \eref{2.3}--\eref{2.4} is 
\begin{equation}
\begin{aligned}
A(\mU):=\frac{\partial\mF}{\partial\mU}(\bm U)=\begin{pmatrix}0&1&0\\[0.5ex]\dfrac{\gamma-3}{2}u^2&(3-\gamma)u&\gamma-1\\[1.5ex]
-\dfrac{\gamma uE}{\rho}+(\gamma-1)u^3&H-(\gamma-1)u^2&\gamma u\end{pmatrix},
\end{aligned}
\label{starm1}
\end{equation}
where $H:=(E+p)/\rho$ is the total specific enthalpy. The eigenvalues of $A$ are $\lambda_1(A)=u-c$, $\lambda_2(A)=u$, and
$\lambda_3(A)=u+c$, where $c:=\sqrt{\gamma p/\rho}$ is the speed of sound.

\subsection{One-Dimensional Second-Order Adaptive Artificial Anti-Diffusion Scheme}\label{sec21}
Let the computational domain be partitioned into uniform cells $I_j:=[x_\jmh,x_\jph]$ of size $x_\jph-x_\jmh\equiv\dx$, centered at
$x_j=\big(x_\jmh+x_\jph\big)/2$. Assume that the cell averages,
$$
\xbar\mU_j(t)\approx\frac{1}{\dx}\int\limits_{I_j}\mU(x,t)\,{\rm d}x,
$$
are available at a certain time $t\ge0$. In what follows, for the sake of brevity, we suppress the time dependence of indexed quantities.

In the semi-discrete framework, the cell averages are evolved by numerically solving the system of ODEs
\begin{equation}
\frac{{\rm d}\xbar{\mU}_j}{{\rm d}t}=-\frac{\bm{{\cal F}}^{\rm FV}_\jph-\bm{{\cal F}}^{\rm FV}_\jmh}{\dx},
\label{star0}
\end{equation}
where $\bm{{\cal F}}^{\rm FV}_\jph$ are FV numerical fluxes.

Let us assume that the numerical fluxes in \eref{star0} are constructed using a nonlinear limiter, which ensures the nonlinear stability of
the resulting FV scheme; for many examples of reliable numerical fluxes see, e.g., the monographs \cite{GR3,Hes,Leveque02,Toro2009} and
references therein. Our goal is to reduce the amount of numerical diffusion present in the numerical flux $\bm{{\cal F}}^{\rm FV}_\jph$ by
adding an AD term to it. The modified semi-discrete scheme will then read as 
\begin{equation*}
\frac{{\rm d}\xbar{\mU}_j}{{\rm d}t}=-\frac{\bm{{\cal F}}^{\rm AD}_\jph-\bm{{\cal F}}^{\rm AD}_\jmh}{\dx},
\end{equation*}
where the numerical flux is given by  
\begin{equation*}
\bm{{\cal F}}^{\rm AD}_\jph=\bm{{\cal F}}^{\rm FV}_\jph-Q_\jph\frac{\,\xbar\mU_{j+1}-\,\xbar\mU_j}{\dx}.
\end{equation*}
Here, the matrix $Q_\jph$ is constructed to primarily add the AD into the linearly degenerate field, which corresponds to the second
eigenvalue $\lambda_2=u$ of the Jacobian \eref{starm1}. To this end, we use the local characteristic decomposition (LCD) based on the
linearization of the Jacobian $A(\bm U)$ about the cell interface point $x=x_\jph$.

We introduce the following local linearization: $\widehat A_\jph=A\big(\widehat{\mU}_\jph\big)$, where the averaged quantities
$\hat{(\cdot)}$ are obtained by averaging the primitive variables, which results in
\begin{equation}
\begin{aligned}
&\hat\rho=\frac{\xbar\rho_j+\xbar\rho_{j+1}}{2},\quad\hat u=\frac{u_j+u_{j+1}}{2},\quad\hat p=\frac{p_j+p_{j+1}}{2},\quad 
\hat E=\frac{\hat p}{\gamma-1}+\hf\hat\rho\hat u^2,\\
&\hat H=\frac{\hat E+\hat p}{\hat\rho},\quad\hat\phi=2\hat H-\hat u^2,\quad\hat c=\sqrt{\gamma\hat p/\hat\rho},
\end{aligned}
\label{star1}
\end{equation}
where the point values of $u$ and $p$ at $x=x_j$ are given by
$$
u_j=\frac{(\xbar{\rho u})_j}{\xbar\rho_j}\quad\mbox{and}\quad p_j=(\gamma-1)\left[\xbar E_j-\hf\rho_ju_j^2\right].
$$
Notice that all of the $\hat{(\cdot)}$ quantities have to have a subscript index $\hat{(\cdot)}_\jph$, but we have omitted it for the sake
of brevity for all of the quantities except for $\widehat A_\jph$. We also stress that the averages in \eref{star1} can be computed in
several alternative ways, for instance, using the Roe averages \cite{Roe81}; see also \cite{GR3,Hes,Leveque02,Toro2009}.

Equipped with $\widehat A_\jph$, we add the AD into the linearly degenerate field by taking a diagonal AD matrix in the local
characteristic space,
\begin{equation}
\Lambda_\jph:=\begin{pmatrix}0&0&0\\0&\texttt C_\jph&0\\0&0&0\end{pmatrix},
\label{star4}
\end{equation}
where $\texttt C_\jph>0$ is yet to be chosen, and then use $Q_\jph=-R_\jph\Lambda_\jph\big(R_\jph\big)^{-1}$. Here, the matrix $R_\jph$
is composed of the eigenvectors of $\widehat A_\jph$ so that $\big(R_\jph\big)^{-1}\widehat A_\jph R_\jph$ is diagonal, and the matrices
$R_\jph$ and $\big(R_\jph\big)^{-1}$ are given by
\begin{equation}
R_\jph=\begin{pmatrix}1&1&1\\\hat u-\hat c&\hat u&\hat u+\hat c\\[0.5ex]\hat H-\hat u\hat c&\dfrac{\hat u^2}{2}&\hat H+\hat u\hat c
\end{pmatrix},\quad\big(R_\jph\big)^{-1}=\dfrac{1}{\hat\phi}\begin{pmatrix}
\dfrac{\hat u^2}{2}+\dfrac{\hat u\hat\phi}{2\hat c}&-\hat u-\dfrac{\hat\phi}{2\hat c}&1\\[1.5ex]2\hat\phi-{2\hat H}&{2\hat u}&-2\\[0.3ex]
\dfrac{\hat u^2}{2}-\dfrac{\hat u\hat\phi}{2\hat c}&-\hat u+\dfrac{\hat\phi}{2\hat c}&1\end{pmatrix}.
\label{star3}
\end{equation}

The values $\texttt C_\jph$ in \eref{star4} are selected adaptively according to the local smoothness of the computed solution. To this end,
we follow \cite{CHK2026} and evaluate the modified minmod-based SIs. First, we compute 
$$
s^\rho_j=\frac{1}{\max\{\,\xbar\rho_{j-1},\,\xbar\rho_j,\,\xbar\rho_{j+1}\}}\,
{\rm minmod}\left(\,\xbar\rho_{j+1}-\,\xbar\rho_j,\,\xbar\rho_j-\,\xbar\rho_{j-1}\right),
$$
and then check whether
\begin{equation}
|s^\rho_j|>\max\big\{|s^\rho_{j-1}|,|s^\rho_{j+1}|\big\}+\eps_0,
\label{2.5}
\end{equation}
where $\eps_0$ is a positive constant taken to be $\eps_0=0.002$ as in \cite{CHK2026}. If \eref{2.5} is satisfied, then we mark the cells
$I_{j-1}$, $I_j$, and $I_{j+1}$ as ``rough''.

For those cells, at which \eref{2.5} is satisfied, we also compute
$$
s^p_j=\frac{1}{\max\{p_{j-1},p_j,p_{j+1}\}}\,{\rm minmod}\left(p_{j+1}-p_j, p_j-p_{j-1}\right),
$$
and check whether
\begin{equation}
|s^p_j|<\max\big\{|s^p_{j-1}|,|s^p_{j+1}|\big\},
\label{2.6}
\end{equation}
and if \eref{2.6} is satisfied as well, then we mark the cells $I_{j-1}$, $I_j$, and $I_{j+1}$ as ``rough contact'' cells.

Finally, the cells that are not marked as either ``rough'' or ``rough contact'' are marked as ``smooth''. 

After all of the cells have been identified as either ``rough'', ``rough contact'', or ``smooth'', we select $\texttt C_\jph$ according to
Algorithm 1.
\begin{algorithm}[ht!]
\caption{Assignment of $\texttt C_\jph$}
\begin{algorithmic}[1]
\For{each cell interface $x=x_\jph$}
\If{either $I_j$ or $I_{j+1}$ is a ``rough contact'' cell}
\State set $\texttt C_\jph=\texttt C\,\dx$
\Else
\State set $\texttt C_\jph=\texttt C\,(\dx)^2$
\EndIf
\EndFor
\end{algorithmic}
\end{algorithm}
\begin{remark}
The parameter $\texttt C$ is a tunable constant that should be selected for each problem at hand. In practice, we tune $\texttt C$ on a
coarse mesh and then use the same value on finer meshes. An example of this tuning procedure is provided in Example 4 in \S\ref{sec41}.
\end{remark}
\begin{remark}
The proposed AAAD method can be extended to other hyperbolic systems of conservation laws possessing one or more linearly degenerate
characteristic fields. In such cases, the AD correction can be introduced to sharpen linearly degenerate waves without affecting the
genuinely nonlinear fields.
\end{remark}

\subsection{One-Dimensional Fifth-Order Adaptive Artificial Anti-Diffusion Scheme}\label{sec22}
We now construct a fifth-order AAAD scheme within the A-WENO framework, which we briefly review in \S\ref{sec221} and then introduce the
fifth-order adaptive algorithm in \S\ref{sec222}.

\subsubsection{One-Dimensional Fifth-Order A-WENO Schemes}\label{sec221}
Following \cite{JSZ}, the point values $\mU_j\approx\bm U(x_j,t)$ are evolved in time by numerically solving
\begin{equation}
\frac{{\rm d}\mU_j}{{\rm d}t}=-\frac{\bm{{\cal H}}_\jph-\bm{{\cal H}}_\jmh}{\dx},
\label{star2}
\end{equation}
where the fifth-order numerical flux is
\begin{equation*}
\bm{{\cal H}}_\jph=\bm{{\cal F}}^{\rm FV}_\jph-\frac{1}{24}(\dx)^2(\mF_{xx})_\jph+\frac{7}{5760}(\dx)^4(\mF_{xxxx})_\jph.
\end{equation*}
Here, $\bm{{\cal F}}^{\rm FV}_\jph$ is a FV numerical flux, and $(\mF_{xx})_\jph$ and $(\mF_{xxxx})_\jph$ are higher-order correction terms
computed using fourth- and second-order central differences, respectively:
$$
\begin{aligned}
&(\mF_{xx})_\jph=\frac{1}{48(\dx)^2}\Big[-5\mF_{j-2}+39\mF_{j-1}-34\mF_j-34\mF_{j+1}+39\mF_{j+2}-5\mF_{j+3}\Big],\\
&(\mF_{xxxx})_\jph=\frac{1}{2(\dx)^4}\Big[\mF_{j-2}-3\mF_{j-1}+2\mF_j+2\mF_{j+1}-3\mF_{j+2}+\mF_{j+3}\Big],
\end{aligned}
$$
where $\mF_j:=\mF(\mU_j)$. 

In general, the FV fluxes are to be computed using the one-sided point values $\mU^\pm_\jph$, which must be reconstructed with at least
fifth order of accuracy. We employ the fifth-order WENO-Z interpolation \cite{JSZ,Liu17} applied to the local characteristic variables; see
Appendix \ref{appc} for details.

\subsubsection{One-Dimensional Fifth-Order Adaptive Artificial Anti-Diffusion Algorithm}\label{sec222}
The fifth-order adaptive scheme is obtained by replacing the numerical fluxes $\bm{{\cal H}}_{j\pm\hf}$ in \eref{star2} with
$\bm{{\cal H}}^{\rm AD}_{j\pm\hf}$:
\begin{equation*}
\frac{{\rm d}\mU_j}{{\rm d}t}=-\frac{\bm{{\cal H}}^{\rm AD}_\jph-\bm{{\cal H}}^{\rm AD}_\jmh}{\dx},
\end{equation*}
where the adaptive numerical flux is
\begin{equation*}
\bm{{\cal H}}^{\rm AD}_\jph=\bm{{\cal H}}_\jph-Q_\jph\frac{\mU_{j+1}-\mU_j}{\dx},
\end{equation*}
and the matrix $Q_\jph$ is defined as in the second-order AAAD scheme (see \S\ref{sec21}) but with the values of $\texttt C_\jph$
determined according to Algorithm 2.
\begin{algorithm}[ht!]
\caption{Assignment of $\texttt C_\jph$}
\begin{algorithmic}[1]
\For{each cell interface $x=x_\jph$}
\If{either $I_j$ or $I_{j+1}$ is a ``rough contact'' cell}
\State set $\texttt C_\jph=\texttt C\,\dx$
\ElsIf{either $I_j$ or $I_{j+1}$ is a ``rough'' cell}
\State set $\texttt C_\jph=\texttt C\,(\dx)^2$
\Else
\State set $\texttt C_\jph=\texttt C\,(\dx)^5$
\EndIf
\EndFor
\end{algorithmic}
\end{algorithm}

As in the second-order case, this choice ensures that the added AD term does not affect the designed order of accuracy in ``smooth''
regions, while enhancing the resolution of linearly degenerate waves.

\section{Two-Dimensional Adaptive Artificial Anti-Diffusion Schemes}\label{sec3}
In this section, we introduce the second- and fifth-order AAAD methods for the 2-D Euler equations of gas dynamics, which read as \eref{1.2}
with
\begin{equation}
\mU=(\rho,\rho u,\rho v,E)^\top,\quad\mF=\big(\rho u,\rho u^2+p,\rho uv,u(E+p)\big)^\top,\quad
\mG=\big(\rho v,\rho uv,\rho v^2+p,v(E+p)\big)^\top.
\label{3.1f}
\end{equation}
Here, $u$ and $v$ are the $x$- and $y$-directional velocities, respectively, and the pressure is given by the ideal-gas EOS
\begin{equation}
p=(\gamma-1)\Big[E-\frac{\rho}{2}(u^2+v^2)\Big].
\label{3.2f}
\end{equation}

An explicit form of the Jacobians $A(\mU):=\frac{\partial\mF}{\partial\mU}$ and $B(\mU):=\frac{\partial\mG}{\partial\mU}$ of the system
\eref{1.2}, \eref{3.1f}--\eref{3.2f} can be found in, e.g., \cite{Toro2009}. Their eigenvalues are $\lambda_1(A)=u-c$,
$\lambda_2(A)=\lambda_3(A)=u$, $\lambda_4(A)=u+c$, and $\lambda_1(B)=v-c$, $\lambda_2(B)=\lambda_3(B)=v$, $\lambda_4(B)=v+c$, respectively.

\subsection{Two-Dimensional Second-Order Adaptive Artificial Anti-Diffusion Scheme}\label{sec31}
Let the computational domain be partitioned into uniform cells $I_{j,k}:=[x_\jmh,x_\jph]\times[y_\kmh,y_\kph]$, centered at
$(x_j,y_k)=\big(\big(x_\jmh+x_\jph\big)/2,\big(y_\kmh+y_\kph\big)/2\big)$, with $x_\jph-x_\jmh\equiv\dx$ and $y_\kph-y_\kmh\equiv\dy$.

Assume that the computed cell averages,
$$
\xbar\mU_{j,k}\approx\frac{1}{\dx\dy}\int\limits_{I_{j,k}}\mU(x,y,t)\,{\rm d}x\,{\rm d}y,
$$
are available at a certain time level $t\ge0$. They are evolved by numerically solving the following semi-discrete system:
\begin{equation*}
\frac{{\rm d}\xbar{\mU}_{j,k}}{{\rm d}t}=-\frac{\bm{{\cal F}}^{\rm AD}_{\jph,k}-\bm{{\cal F}}^{\rm AD}_{\jmh,k}}{\dx}-
\frac{\bm{{\cal G}}^{\rm AD}_{j,\kph}-\bm{{\cal G}}^{\rm AD}_{j,\kmh}}{\dy},
\end{equation*}
where the adaptive numerical fluxes are defined by
\begin{equation*}
\bm{{\cal F}}^{\rm AD}_{\jph,k}=\bm{{\cal F}}^{\rm FV}_{\jph,k}-Q_{\jph,k}\frac{\,\xbar\mU_{j+1,k}-\,\xbar\mU_{j,k}}{\dx},\qquad 
\bm{{\cal G}}^{\rm AD}_{j,\kph}=\bm{{\cal G}}^{\rm FV}_{j,\kph}-Q_{j,\kph}\frac{\,\xbar\mU_{j,k+1}-\,\xbar\mU_{j,k}}{\dy}.
\end{equation*}
Here, $\bm{{\cal F}}^{\rm FV}_{\jph,k}$ and $\bm{{\cal G}}^{\rm FV}_{j,\kph}$ are FV numerical fluxes. As in the 1-D case, we assume that
these numerical fluxes are constructed using a nonlinear limiter, which ensures the nonlinear stability of the underlying FV scheme.

To construct the matrices $Q_{\jph,k}$ and $Q_{j,\kph}$, we first introduce the linearizations
$\widehat A_{\jph,k}=A\big(\widehat\mU_{\jph,k}\big)$ and $\widehat B_{j,\kph}=B\big(\widehat\mU_{j,\kph}\big)$, where
$\widehat\mU_{\jph,k}$ is a proper average of $\,\xbar{\bm U}_{j,k}$ and $\,\xbar{\bm U}_{j+1,k}$ and $\widehat\mU_{j,\kph}$ is a proper
average of $\,\xbar{\bm U}_{j,k}$ and $\,\xbar{\bm U}_{j,k+1}$. We then compute the matrices $R_{\jph,k}$ and $R_{j,\kph}$, for which
$\big(R_{\jph,k}\big)^{-1}\widehat A_{\jph,k}R_{\jph,k}$ and $\big(R_{j,\kph}\big)^{-1}\widehat B_{j,\kph}R_{j,\kph}$ are diagonal.
\begin{remark}
The $x$-dimensional averages $\widehat\mU_{\jph,k}$ as well as the matrices $R_{\jph,k}$ and $\big(R_{\jph,k}\big)^{-1}$ can be found in
\cite[Appendix C]{CCHKL_22}. The $y$-dimensional averages $\widehat\mU_{j,\kph}$ and the matrices $R_{j,\kph}$ and
$\big(R_{j,\kph}\big)^{-1}$ can be obtained in a similar way.
\end{remark}

As in the 1-D case, we add the AD terms to the linearly degenerate fields only. Specifically, we take
$Q_{\jph,k}=-R_{\jph,k}\Lambda_{\jph,k}\big(R_{\jph,k}\big)^{-1}$ and $Q_{j,\kph}=-R_{j,\kph}\Lambda_{j,\kph}\big(R_{j,\kph}\big)^{-1}$,
where 
$$
\Lambda_{\jph,k}:=\begin{pmatrix}0&0&0&0\\0&\texttt C_{\jph,k}&0&0\\0&0&\texttt C_{\jph,k}&0\\0&0&0&0\end{pmatrix},\quad
\Lambda_{j,\kph}:=\begin{pmatrix}0&0&0&0\\0&\texttt C_{j,\kph}&0&0\\0&0&\texttt C_{j,\kph}&0\\0&0&0&0\end{pmatrix}.
$$

The values of $\texttt C_{\jph,k}$ and $\texttt C_{j,\kph}$ are selected adaptively according to the local smoothness of the computed
solution. To this end, we follow \cite{CHK2026} and evaluate modified minmod-based SIs for $\rho$ and $p$ in a dimension-by-dimension
manner. We now show how to check the local smoothness of the computed solution and to determine the values of $\texttt C_{\jph,k}$ in the
$x$-direction. (The values of $\texttt C_{j,\kph}$ can be obtained similarly, and we omit the details for the sake of brevity.)

We first compute the $x$-directional SIs:
$$
s^{x,\rho}_{j,k}=
\frac{1}{\max\{\,\xbar\rho_{j-1,k},\,\xbar\rho_{j,k},\,\xbar\rho_{j+1,k}\}}\,
{\rm minmod}\left(\,\xbar\rho_{j+1,k}-\,\xbar\rho_{j,k},\,\xbar\rho_{j,k}-\,\xbar\rho_{j-1,k}\right),
$$
and then check whether
\begin{equation}
|s^{x,\rho}_{j,k}|>\max\big\{|s^{x,\rho}_{j-1,k}|,|s^{x,\rho}_{j+1,k}|\big\}+\eps_0,
\label{4.5}
\end{equation}
where $\eps_0$ is a positive constant taken to be $\eps_0=0.002$ as in the 1-D case. If \eref{4.5} is satisfied, then we mark the cells
$I_{j-1,k}$, $I_{j,k}$, and $I_{j+1,k}$ as ``rough''.

For those cells, at which \eref{4.5} is satisfied, we also compute the SIs:
$$
s^{x,p}_{j,k}=\frac{1}{\max\{p_{j-1,k},p_{j,k},p_{j+1,k}\}}\,{\rm minmod}\left(p_{j+1,k}-p_{j,k},p_{j,k}-p_{j-1,k}\right),
$$
and check whether 
\begin{equation}
|s^{x,p}_{j,k}|<\max\big\{|s^{x,p}_{j-1,k}|,|s^{x,p}_{j+1,k}|\big\},
\label{4.6}
\end{equation}
and if \eref{4.6} is satisfied as well, then we mark the cells $I_{j-1,k}$, $I_{j,k}$, and $I_{j+1,k}$ as ``rough contact'' cells.

Finally, the cells that are not marked as either ``rough'' or ``rough contact'' are marked as ``smooth''. 

After all of the cells have been identified as either ``rough'', ``rough contact'', or ``smooth'', we select $\texttt C_{\jph,k}$ according
to Algorithm 3.
\setcounter{algorithm}{2}
\begin{algorithm}[ht!]
\caption{Assignment of $\texttt C_{\jph,k}$}
\begin{algorithmic}[1]
\For{each $(x_\jph,y_k)$}
\If{either $I_{j,k}$ or $I_{j+1,k}$ is a ``rough contact'' cell}
\State set $\texttt C_{\jph,k}=\texttt C\,\dx$
\Else
\State set $\texttt C_{\jph,k}=\texttt C\,(\dx)^2$
\EndIf
\EndFor
\end{algorithmic}
\end{algorithm}
\begin{remark}
As in the 1-D case, $\texttt C$ is a positive tunable constant to be selected for each problem at hand. In practice, one may tune
$\texttt C$ on a coarse mesh and then use the same value on finer meshes, as demonstrated in Example 9 in \S\ref{sec42}. 
\end{remark}

\subsection{Two-Dimensional Fifth-Order Adaptive Artificial Anti-Diffusion Scheme}
We now extend the 2-D second-order AAAD method to the fifth order of accuracy via the framework of the A-WENO scheme. The point values
$\mU_{j,k}\approx\mU(x_j,y_k,t)$ are evolved in time by numerically solving
\begin{equation*}
\frac{{\rm d}\mU_{j,k}}{{\rm d}t} =-\frac{\bm{{\cal H}}^{\rm AD}_{\jph,k}-\bm{{\cal H}}^{\rm AD}_{\jmh,k}}{\dx}
-\frac{\bm{{\cal H}}^{\rm AD}_{j,\kph}-\bm{{\cal H}}^{\rm AD}_{j,\kmh}}{\dy},
\end{equation*}
where the numerical fluxes are defined by
\begin{equation*}
\begin{aligned}
\bm{{\cal H}}^{\rm AD}_{\jph,k}&=\bm{{\cal F}}^{\rm FV}_{\jph,k}-\frac{1}{24}(\dx)^2(\mF_{xx})_{\jph,k}+
\frac{7}{5760}(\dx)^4(\mF_{xxxx})_{\jph,k}-Q_{\jph,k}\,\frac{\mU_{j+1,k}-\mU_{j,k}}{\dx},\\
\bm{{\cal H}}^{\rm AD}_{j,\kph}&=\bm{{\cal G}}^{\rm FV}_{j,\kph}-\frac{1}{24}(\dy)^2(\mG_{yy})_{j,\kph}+
\frac{7}{5760}(\dy)^4(\mG_{yyyy})_{j,\kph}-Q_{j,\kph}\,\frac{\mU_{j,k+1}-\mU_{j,k}}{\dy}.
\end{aligned}
\end{equation*}
Here, $\bm{{\cal F}}^{\rm FV}_{\jph,k}$ and $\bm{{\cal G}}^{\rm FV}_{j,\kph}$ are the FV fluxes, and the higher-order correction terms are
computed by fourth- and second-order central differences, respectively:
$$
\begin{aligned}
&(\mF_{xx})_{\jph,k}=\frac{1}{48(\dx)^2}\Big(-5\mF_{j-2,k}+39\mF_{j-1,k}-34\mF_{j,k}-34\mF_{j+1,k}+39\mF_{j+2,k}-5\mF_{j+3,k}\Big),\\
&(\mF_{xxxx})_{\jph,k}=\frac{1}{2(\dx)^4}\Big(\mF_{j-2,k}-3\mF_{j-1,k}+2\mF_{j,k}+2\mF_{j+1,k}-3\mF_{j+2,k}+\mF_{j+3,k}\Big),\\
&(\mG_{yy})_{j,\kph}=\frac{1}{48(\dy)^2}\Big(-5\mG_{j,k-2}+39\mG_{j,k-1}-34\mG_{j,k}-34\mG_{j,k+1}+39\mG_{j,k+2}-5\mG_{j,k+3}\Big),\\
&(\mG_{yyyy})_{j,\kph}=\frac{1}{2(\dy)^4}\Big(\mG_{j,k-2}-3\mG_{j,k-1}+2\mG_{j,k}+2\mG_{j,k+1}-3\mG_{j,k+2}+\mG_{j,k+3}\Big),
\end{aligned}
$$
where $\mF_{j,k}:=\mF(\mU_{j,k})$ and $\mG_{j,k}:=\mG(\mU_{j,k})$.

In general, the FV fluxes are to be computed using the one-sided point values $\mU^\pm_{\jph,k}$ and $\mU^\pm_{j,\kph}$, which must be
reconstructed with at least fifth order of accuracy. To obtain these values, we employ the 1-D fifth-order WENO-Z interpolation (described
in Appendix \ref{appc}) applied to the local characteristic variables in the $x$- and $y$-directions, respectively; we omit the details for
the sake of brevity.

The matrix $Q_{\jph,k}$ is defined as in the second-order AAAD scheme (see \S\ref{sec31}) but with the values of $\texttt C_{\jph,k}$
determined according to Algorithm 4.
\begin{algorithm}[ht!]
\caption{Assignment of $\texttt C_{\jph,k}$}
\begin{algorithmic}[1]
\For{each $(x_\jph,y_k)$}
\If{either $I_{j,k}$ or $I_{j+1,k}$ is a ``rough contact'' cell}
\State set $\texttt C_{\jph,k}=\texttt C\,\dx$
\ElsIf{either $I_{j,k}$ or $I_{j+1,k}$ is a ``rough'' cell}
\State set $\texttt C_{\jph,k}=\texttt C\,(\dx)^2$
\Else
\State set $\texttt C_{\jph,k}=\texttt C\,(\dx)^5$
\EndIf
\EndFor
\end{algorithmic}
\end{algorithm}

The matrix $Q_{j,\kph}$ can be obtained similarly, and we omit the details for the sake of brevity.


\section{Numerical Examples}\label{sec4}
In this section, we test the proposed AAAD schemes on a series of benchmark problems. The AAAD schemes have been implemented using the 1-D
and 2-D CU numerical fluxes, which are briefly described in Appendices \ref{appa} and \ref{appd}, respectively. In the second-order AAAD
schemes, the cell interface point values are reconstructed using the MinMod2 piecewise linear reconstructions; see Appendices \ref{appb} and
\ref{appe} in the 1-D and 2-D cases, respectively.

We compare the performance of the AAAD schemes with that of the corresponding second-order CU and fifth-order A-WENO schemes. For the sake
of brevity, we refer to the second- and fifth-order AAAD schemes as AAAD2 and AAAD5, respectively.

We numerically integrate the semi-discrete systems by the three-stage third-order strong stability preserving Runge-Kutta (SSP RK3) method
(see, e.g., \cite{Gottlieb11,Gottlieb12}) and use the CFL number $0.4$.

\subsection{One-Dimensional Examples}\label{sec41}
We begin with the 1-D Euler equations of gas dynamics, in which we take $\gamma=1.4$.

\subsubsection*{Example 1---1-D Accuracy Test}
In the first example taken from \cite{KKOKC}, we consider the following smooth initial data:
\begin{equation*}
u(x,0)=\sin\Big(\frac{\pi x}{5}+\frac{\pi}{4}\Big),\quad
\rho(x,0)=\bigg[\frac{\gamma-1}{2\sqrt{\gamma}}\Big(u(x,0)+10\Big)\bigg]^{\frac{2}{\gamma-1}},\quad p(x,0)=\rho^\gamma(x,0),
\end{equation*}
subject to the periodic boundary conditions in the computational domain $[0,10]$. We compute the numerical solution until the final time
$t=0.1$ by the studied AAAD2 and AAAD5 schemes (both with the adaptation constant $\texttt C=0.1$) on a sequence of uniform meshes with
$\dx=1/20$, $1/40$, $1/80$, $1/160$, and $1/320$.

We then compute the $L^1$-errors and estimate the experimental convergence rates using the following Runge formulas based on the solutions
computed on three consecutive uniform grids with mesh sizes $\dx$, $2\dx$, and $4\dx$, denoted by $(\cdot)^{\dx}$, $(\cdot)^{2\dx}$, and
$(\cdot)^{4\dx}$, respectively:
\begin{equation*}
{\rm Error}(\dx)\approx\frac{\delta_{12}^2}{|\delta_{12}-\delta_{24}|},\quad
{\rm Rate}(\dx)\approx\log_2\left(\frac{\delta_{24}}{\delta_{12}}\right),
\end{equation*}
where $\delta_{12}:=\|(\cdot)^{\dx}-(\cdot)^{2\dx}\|_{L^1}$ and $\delta_{24}:=\|(\cdot)^{2\dx}-(\cdot)^{4\dx}\|_{L^1}$. The obtained results
are reported in Table \ref{tab1}, where one can clearly observe the expected second- and fifth-order accuracy of the AAAD2 and AAAD5
schemes, respectively.
\begin{table}[ht!]
\centering
\begin{tabular}{|c|cc|cc|}
\hline
\multirow{2}{*}{$\dx$} & \multicolumn{2}{c|}{AAAD2} & \multicolumn{2}{c|}{AAAD5}\\
\cline{2-5}
& Error & Rate & Error & Rate\\
\hline
$1/80$  & 1.11e-04 & 2.05 & 1.24e-08 & 4.80\\
$1/160$ & 2.32e-05 & 2.14 & 2.81e-10 & 5.12\\
$1/320$ & 6.33e-06 & 2.03 & 9.81e-12 & 4.98\\
\hline 
\end{tabular}
\caption{\sf Example 1: The $L^1$-errors of the density $\rho$ and experimental convergence rates for the AAAD schemes.\label{tab1}}
\end{table}

\begin{remark}
We stress that, in order to achieve the fifth order of accuracy for the AAAD5 scheme, we use smaller time steps with $\dt\sim(\dx)^{5/3}$ in
order to balance the spatial and temporal errors.
\end{remark}

\subsubsection*{Example 2---Titarev-Toro Problem}
In the second 1-D example, we consider the shock-entropy wave interaction problem taken from \cite{Toro2005}; see also
\cite{Shu88,Toro2005a}. The initial conditions,
\begin{equation*}
(\rho,u,p)(x,0)=\begin{cases}(1.51695,0.523346,1.805),&x<-4.5,\\(1+0.1\sin(20x),0,1),&x>-4.5,\end{cases}
\end{equation*}
correspond to a forward-facing shock wave of Mach $1.1$ interacting with high-frequency density perturbations. In this example, the
computational domain is $[-5,5]$ and we impose the free boundary conditions.

\begin{figure}[ht!]
\centerline{\includegraphics[trim=0.8cm 0.3cm 0.8cm 0.6cm, clip, width=6cm]{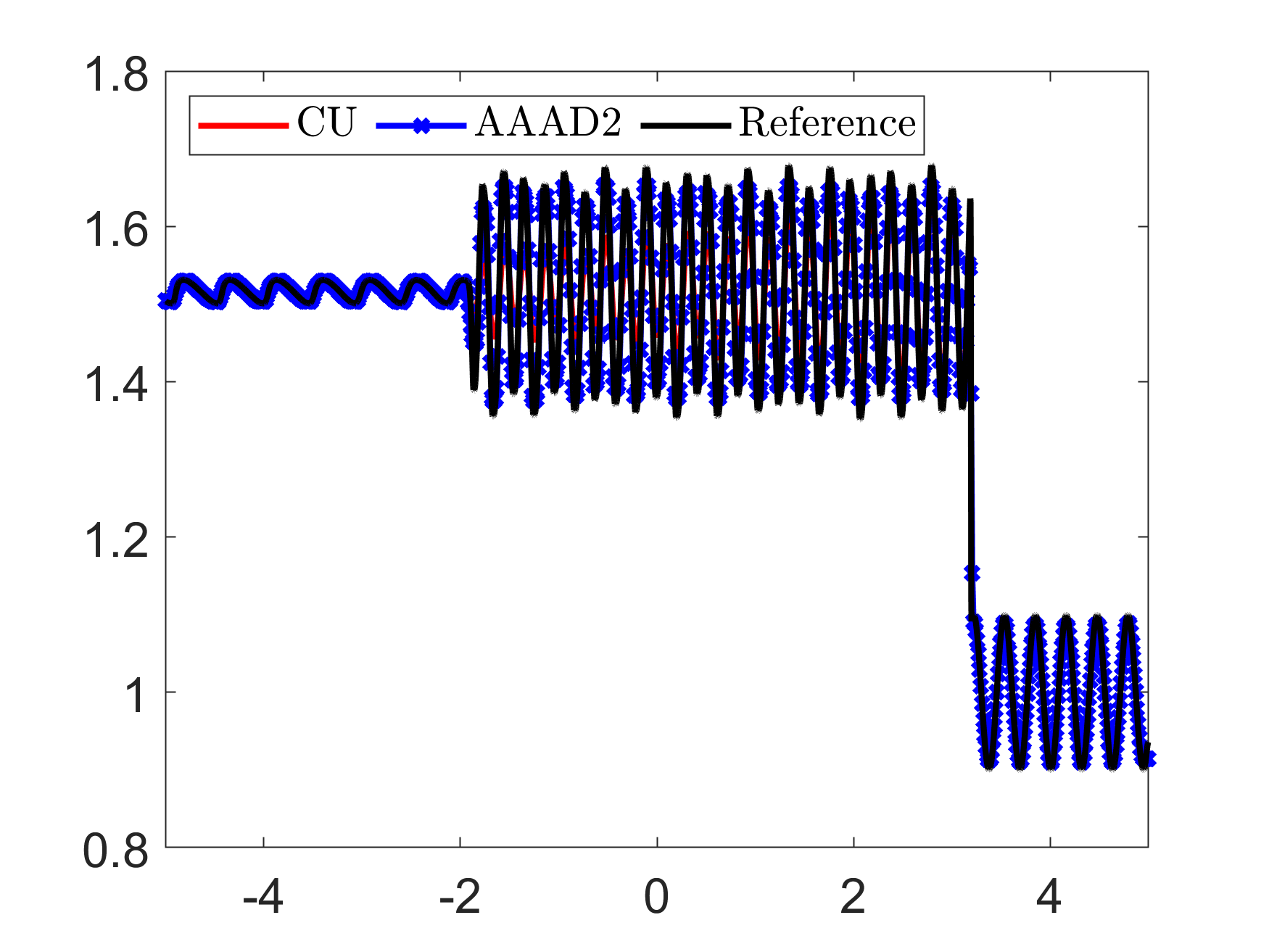}\hspace{1cm}
            \includegraphics[trim=0.8cm 0.3cm 0.8cm 0.6cm, clip, width=6cm]{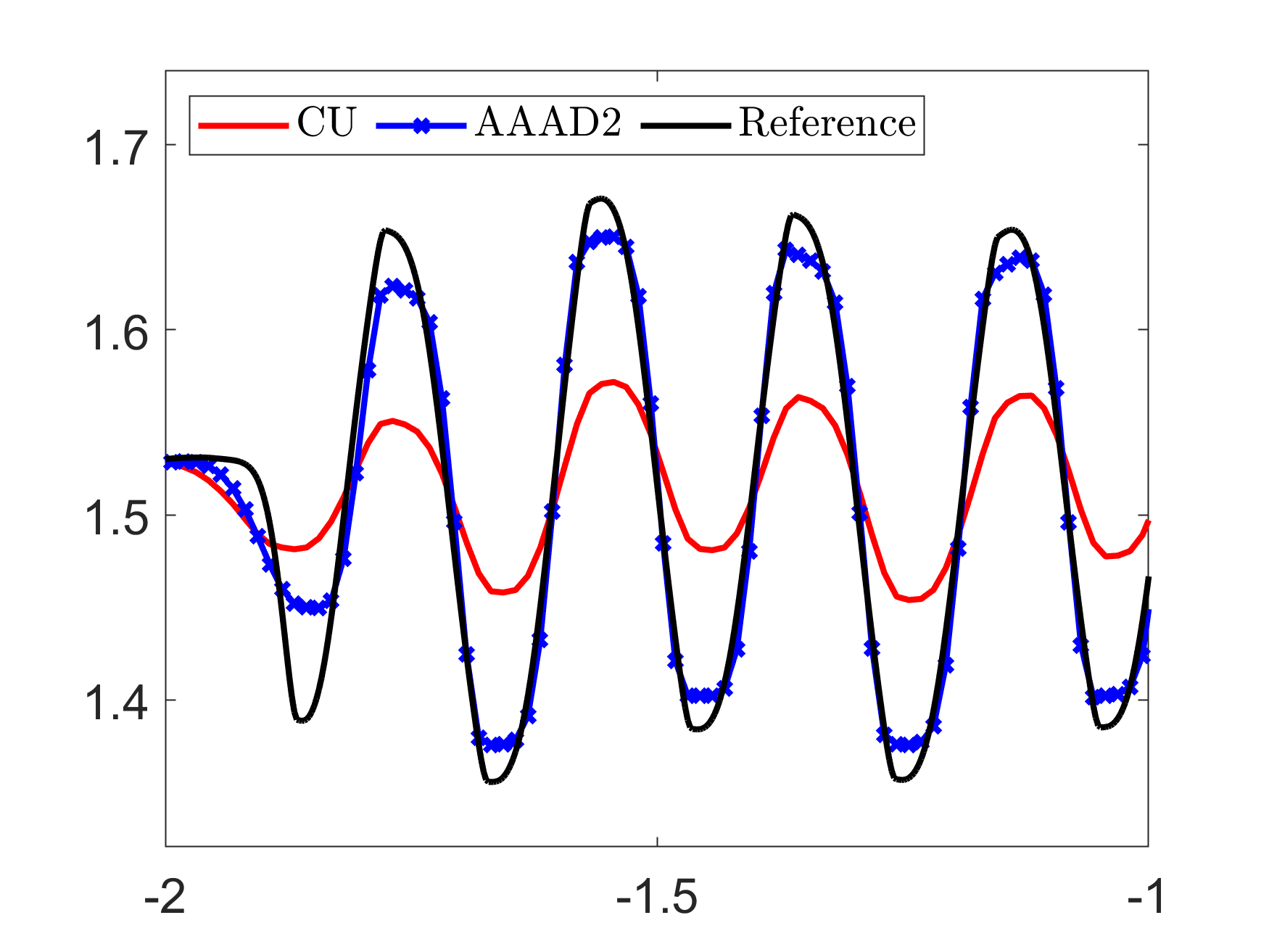}}
\caption{\sf Example 2: Density $\rho$ computed by the CU and AAAD2 schemes on a uniform mesh with $\dx=1/80$ (left) and zoom at
$x\in[-2,-1]$ (right).\label{fig5a}}
\end{figure}

\begin{figure}[ht!]
\centerline{\includegraphics[trim=0.8cm 0.3cm 0.8cm 0.6cm, clip, width=6cm]{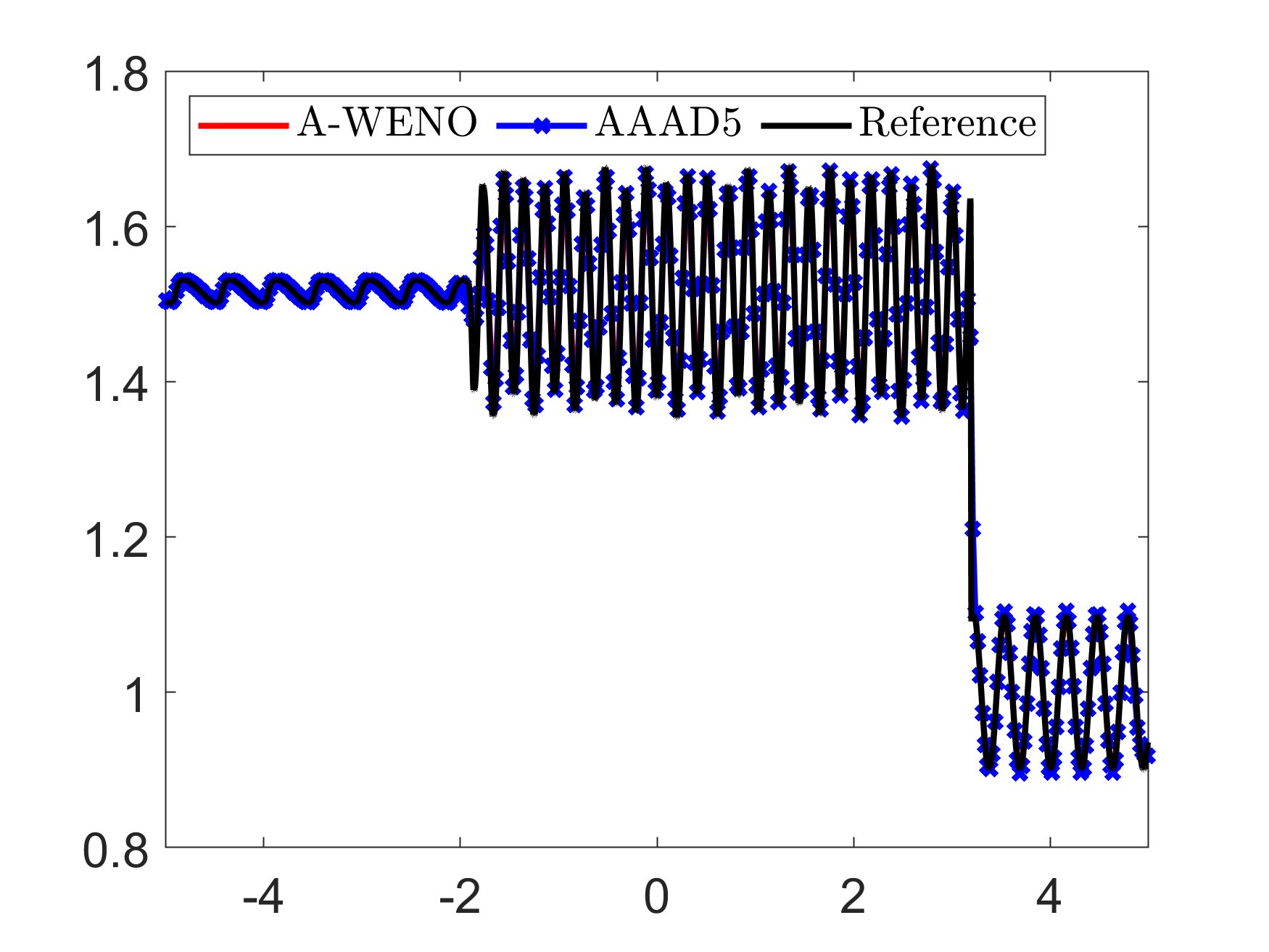}\hspace{1cm}
            \includegraphics[trim=0.8cm 0.3cm 0.8cm 0.6cm, clip, width=6cm]{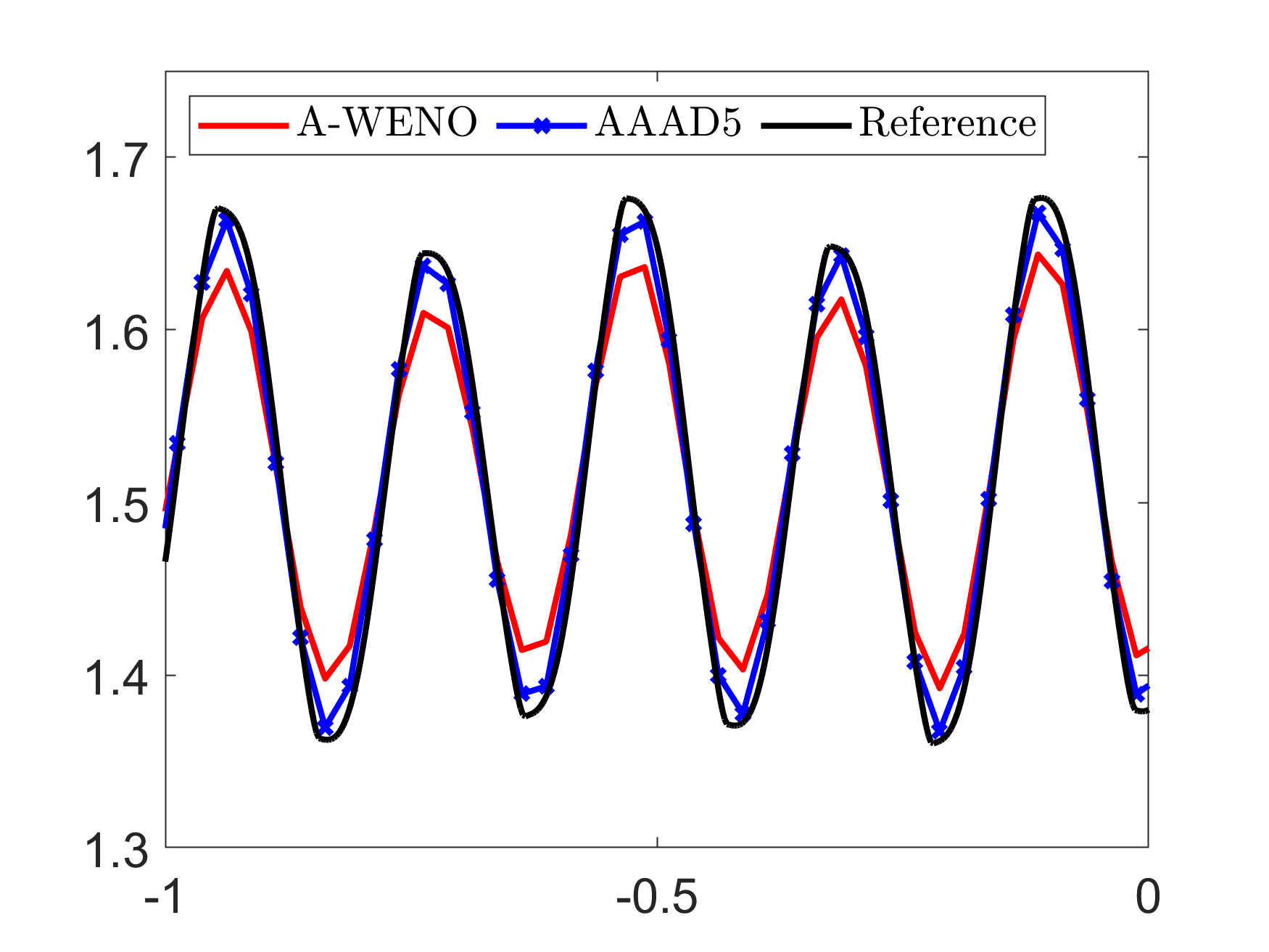}}
\caption{\sf Example 2: Density $\rho$ computed by the A-WENO and AAAD5 schemes on a uniform mesh with $\dx=1/40$ (left), which is coarser
than the mesh used in Figure \ref{fig5a}, and zoom at $x\in[-1,0]$ (right).\label{fig5b}}
\end{figure}

We first compute the numerical solution until the final time $t=5$ by the studied AAAD2 (with the adaptation constant $\texttt C=0.04$) on a
uniform mesh with $\dx=1/80$ and in Figure \ref{fig5a}, we plot the obtained density along with the reference one computed by the CU scheme
on a much finer mesh with $\dx=1/800$. As one can see, the AAAD2 scheme produces much more accurate results than the CU scheme. We then
repeat this computation using the AAAD5 scheme (with the adaptation constant $\texttt C=0.003$), but on a coarser uniform mesh with
$\dx=1/40$. The obtained numerical results, shown in Figure \ref{fig5b}, demonstrate that the fifth-order solution also substantially
improves when the AAAD terms are added. We stress that the fifth-order computations were run on a coarser mesh since finer mesh results were
close to the reference solution. 

\subsubsection*{Example 3---Shock-Density Wave Interaction}
In this example taken from \cite{SO89}, we consider the shock-density wave interaction problem with the initial data,
\begin{equation*}
(\rho,u,p)(x,0)=\begin{cases}\bigg(\dfrac{27}{7},\dfrac{4\sqrt{35}}{9},\dfrac{31}{3}\bigg),&x<-4,\\[0.8ex](1+0.2\sin(5x),0,1),&x>-4,
\end{cases}
\end{equation*}
prescribed in the computational domain $[-5,15]$ with the free boundary conditions.

\begin{figure}[ht!]
\centerline{\includegraphics[trim=0.8cm 0.3cm 0.9cm 0.6cm, clip, width=6cm]{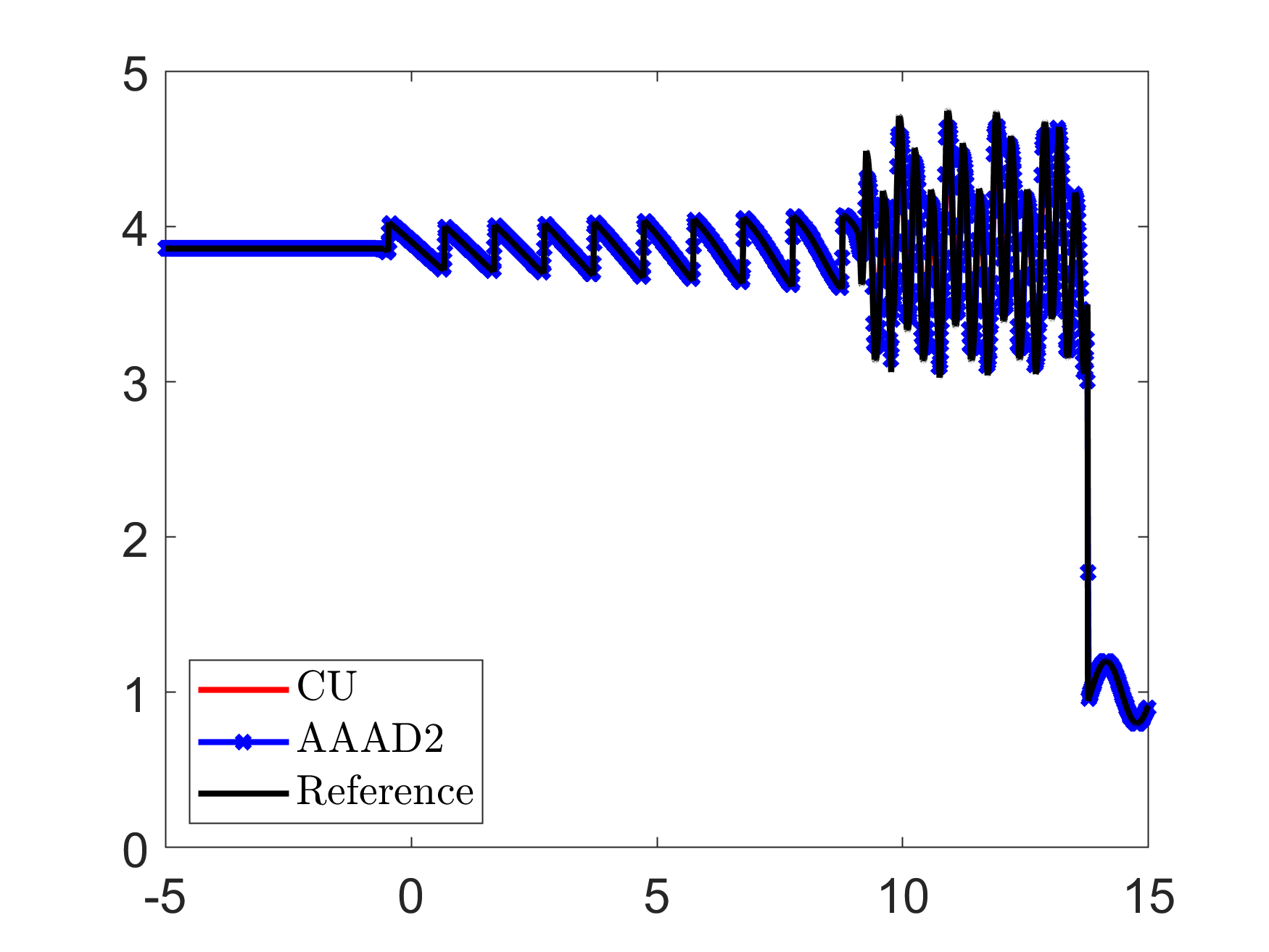}\hspace{1cm}
            \includegraphics[trim=0.8cm 0.3cm 0.9cm 0.6cm, clip, width=6cm]{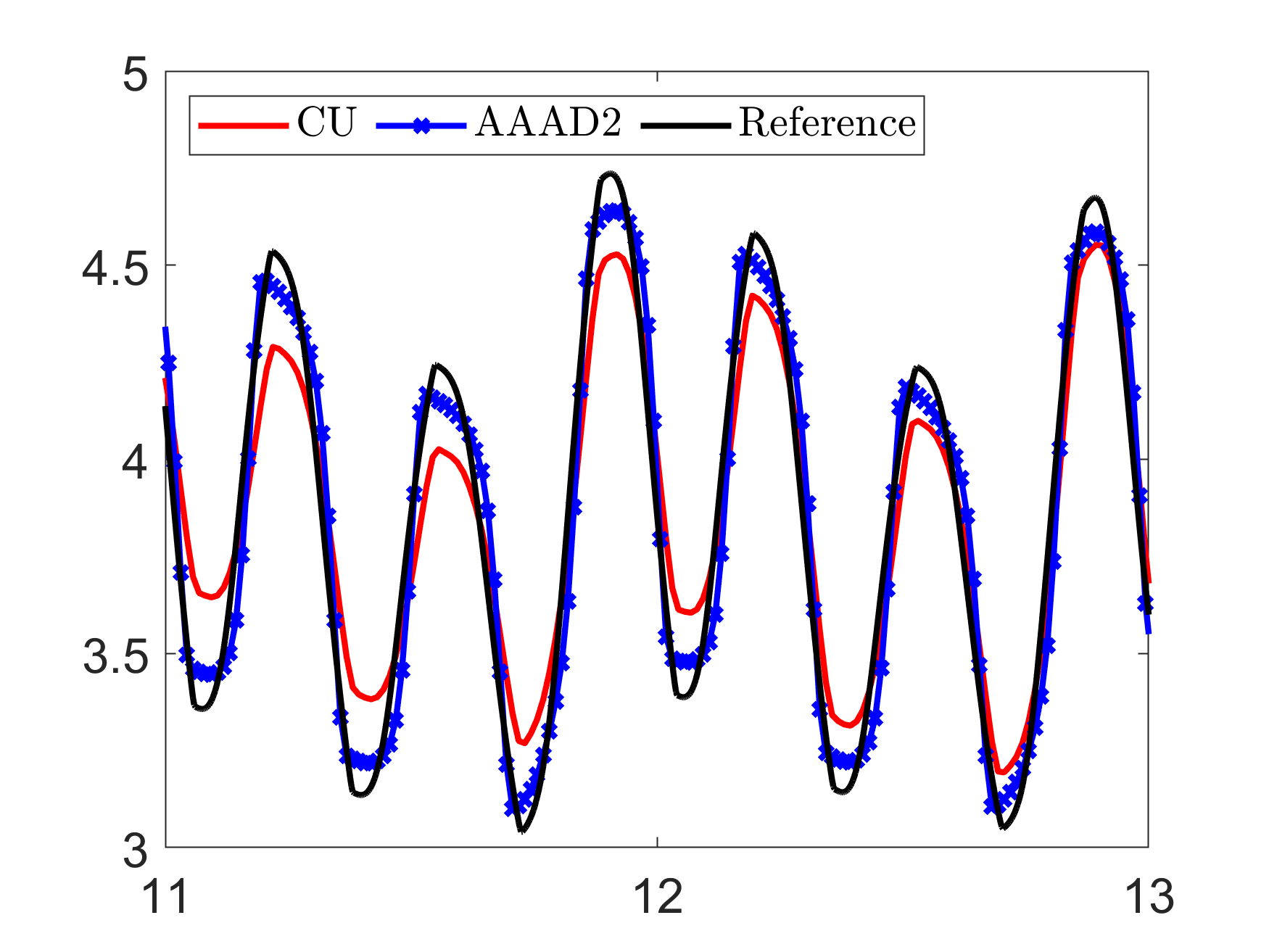}}
\caption{\sf Example 3: Density $\rho$ computed by the CU and AAAD2 schemes on a uniform mesh with $\dx=1/80$ (left) and zoom at
$x\in[11,13]$ (right).\label{fig4a}}
\end{figure}

\begin{figure}[ht!]
\centerline{\includegraphics[trim=0.8cm 0.3cm 0.9cm 0.6cm, clip, width=6cm]{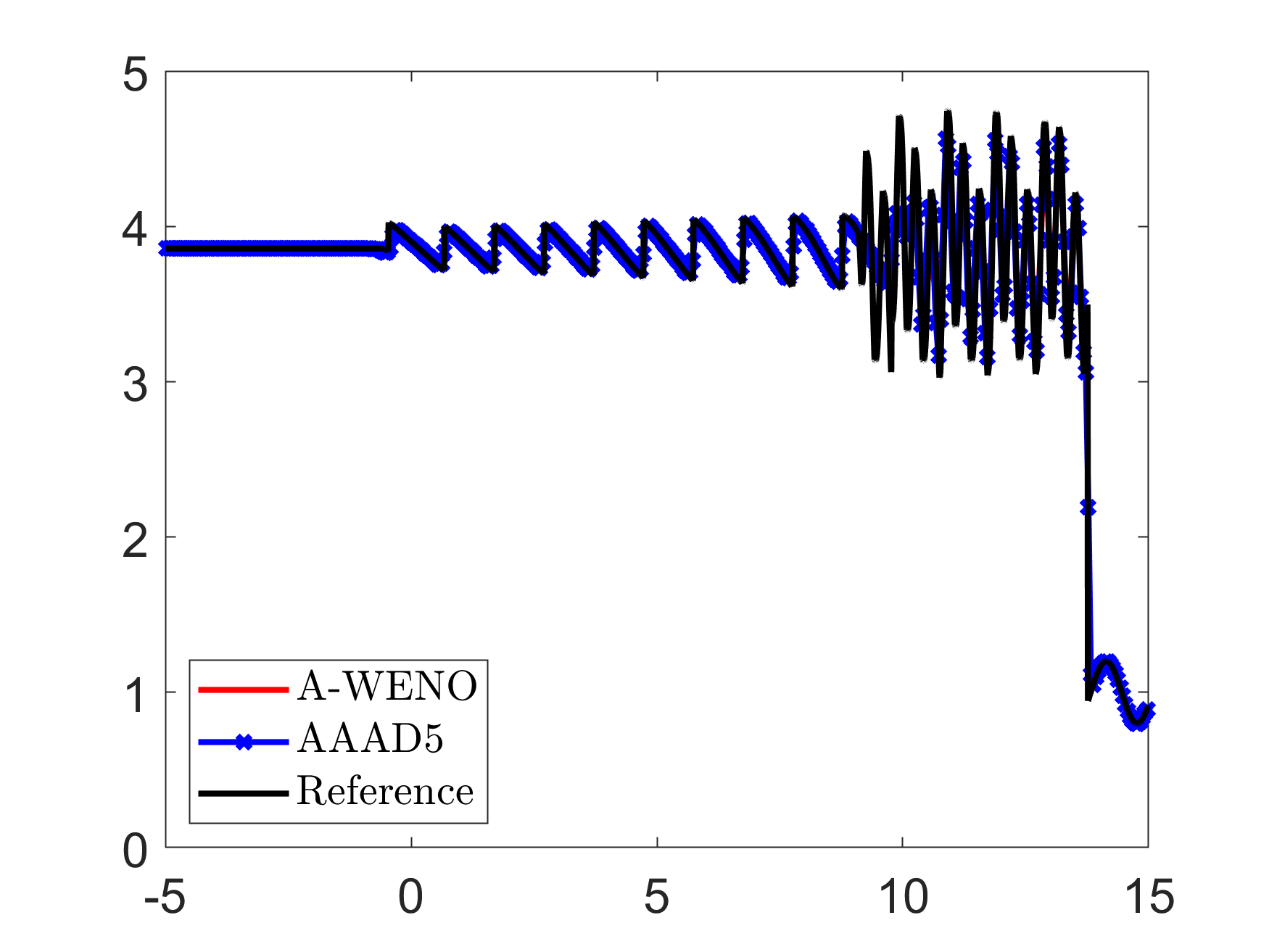}\hspace{1cm}
            \includegraphics[trim=0.8cm 0.3cm 0.9cm 0.6cm, clip, width=6cm]{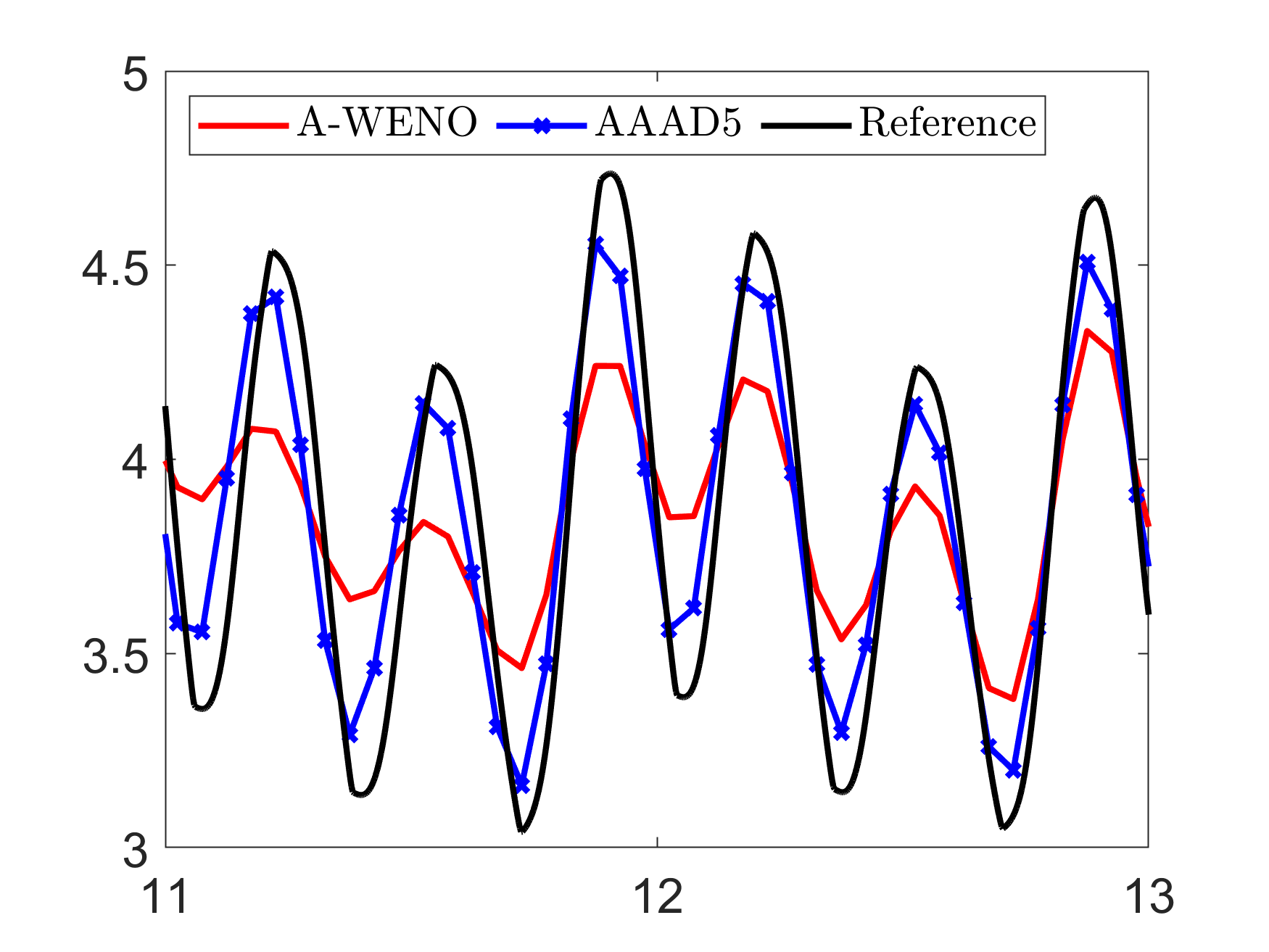}}
\caption{\sf Example 3: Density $\rho$ computed by the A-WENO and AAAD5 schemes on a uniform mesh with $\dx=1/20$ (left), which is much
coarser than the mesh used in Figure \ref{fig4a}, and zoom at $x\in[11,13]$ (right).\label{fig4b}}
\end{figure}

We compute the numerical solutions until the final time $t=5$ by the studied AAAD2 scheme (with the adaptation constant $\texttt C=0.1$) on
a uniform mesh with $\dx=1/80$ and AAAD5 scheme (with the adaptation constant $\texttt C=0.03$) on a coarser uniform mesh with $\dx=1/20$.
We present the results in Figures \ref{fig4a}--\ref{fig4b} together with a reference solution computed by the CU scheme on a much finer mesh
with $\dx=1/400$. One can see that the AAAD2 and AAAD5 schemes yield substantially improved resolution compared with their CU and A-WENO
counterparts.

\subsubsection*{Example 4---Shock-Bubble Interaction}
In the fourth example taken from \cite{KX_22}, we consider the shock-bubble interaction problem with the initial data
\begin{equation*}
(\rho,u,p)(x,0)=\begin{cases}(13.1538,0,1),&\mbox{if }|x|<0.25,\\(1.3333,-0.3535,1.5),&\mbox{if }x>0.75,\\(1,0,1),&\mbox{otherwise},
\end{cases}
\end{equation*}
which correspond to a left-moving shock initially located at $x=0.75$ and a bubble of radius $0.25$ initially located at the origin. We
impose the solid wall boundary conditions on the left and free boundary conditions on the right of the computational domain $[-1,1]$.

\begin{figure}[ht!]
\centerline{\includegraphics[trim=1.0cm 0.3cm 1.3cm 0.8cm, clip, width=6cm]{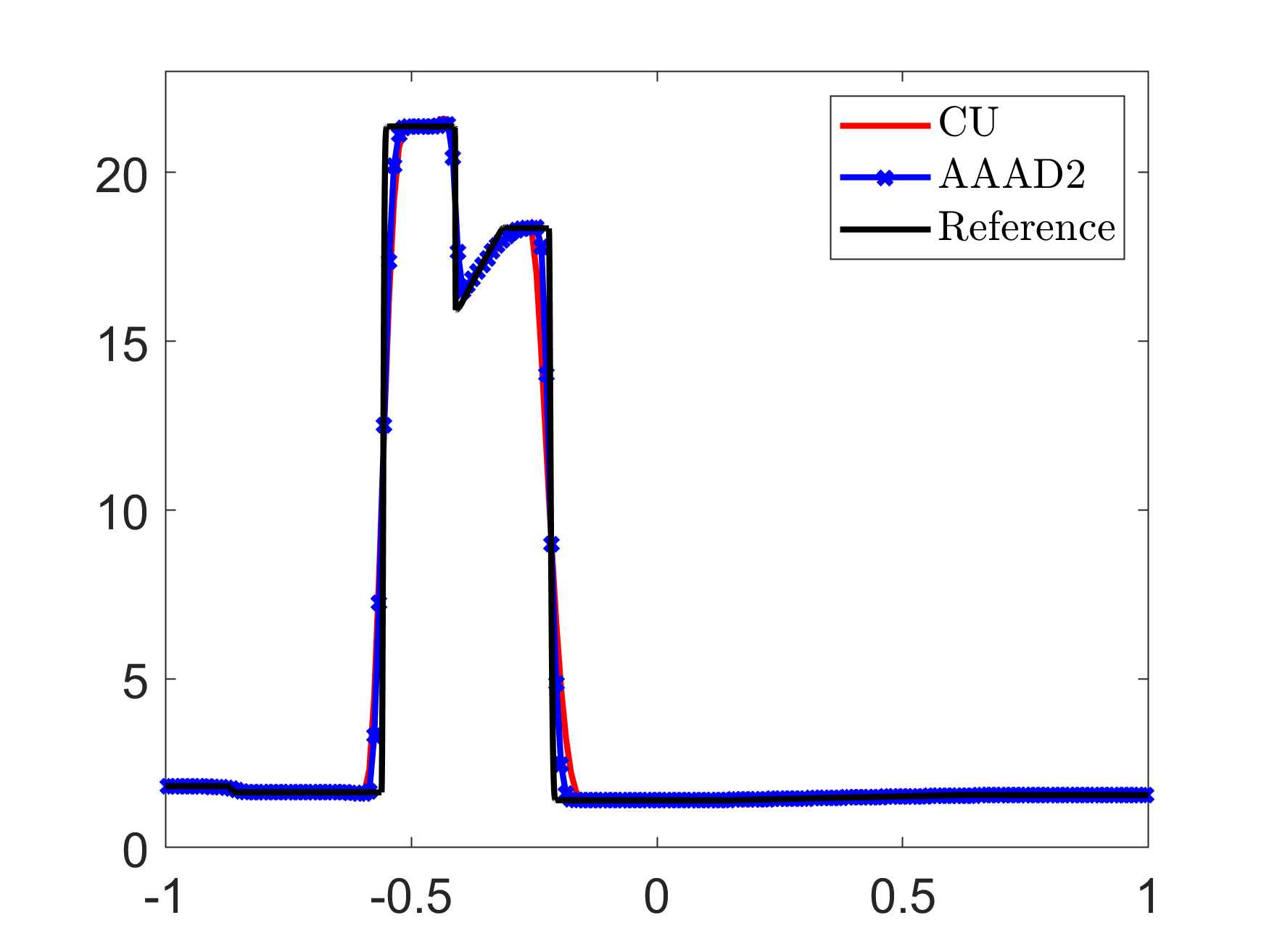}\hspace{1cm}
            \includegraphics[trim=1.0cm 0.3cm 1.3cm 0.8cm, clip, width=6cm]{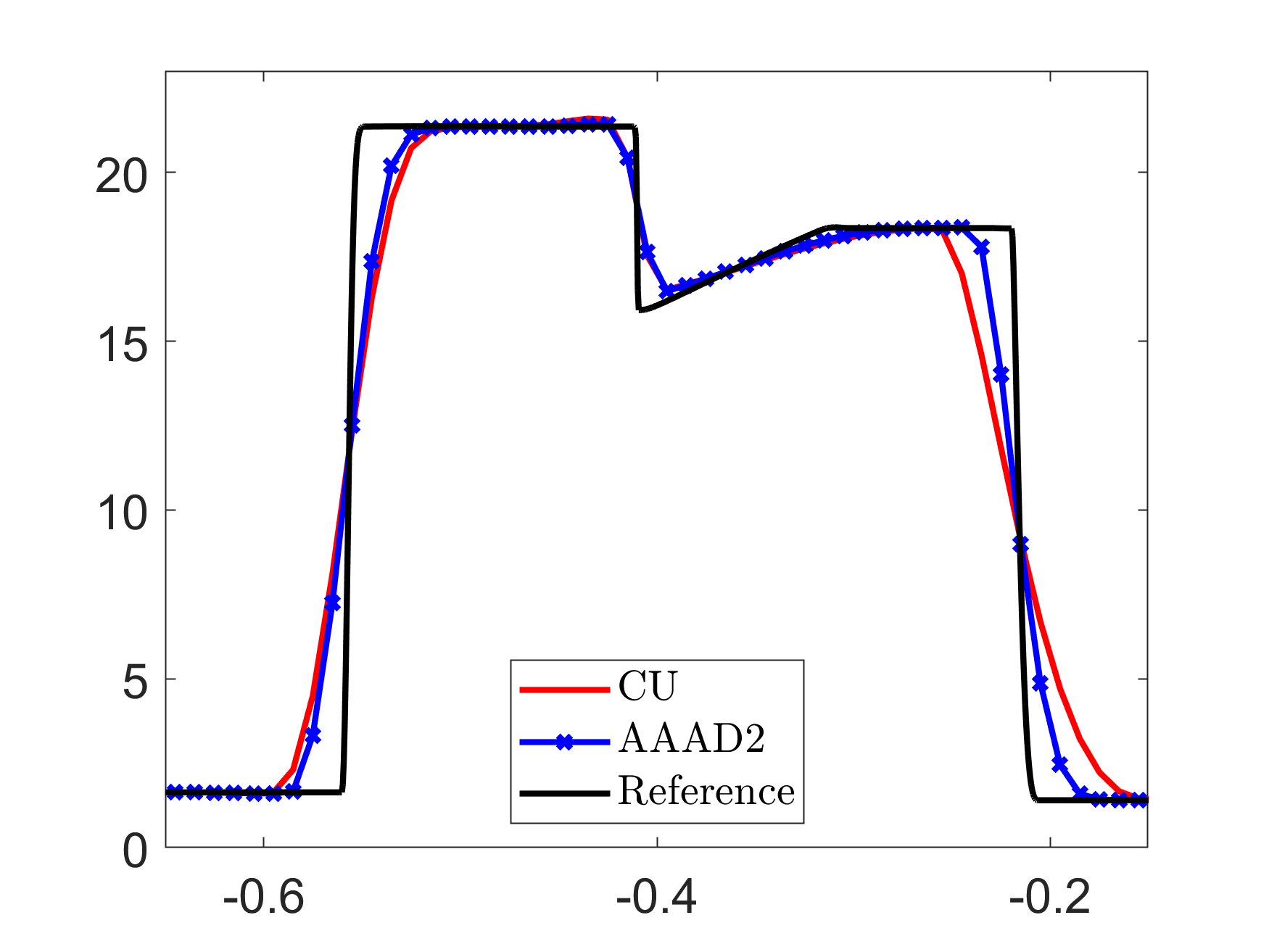}\hspace{1cm}}
\caption{\sf Example 4: Density $\rho$ computed by the CU and AAAD2 schemes (left) and zoom at $x\in[-0.65,-0.15]$ (right).\label{fig3a}}
\end{figure}

\begin{figure}[ht!]
\centerline{\includegraphics[trim=1.0cm 0.3cm 1.3cm 0.8cm, clip, width=6cm]{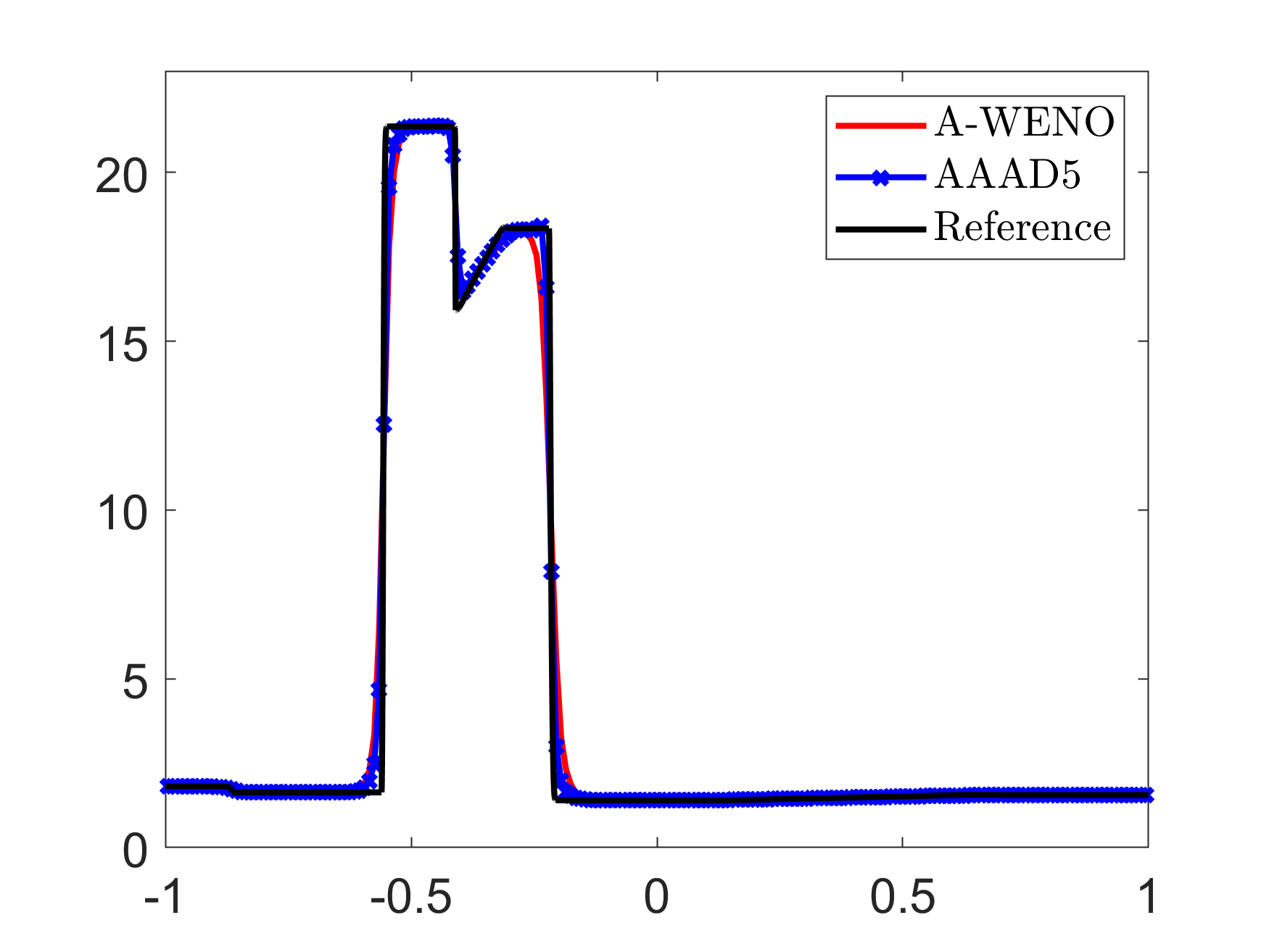}\hspace{1cm}
            \includegraphics[trim=1.0cm 0.3cm 1.3cm 0.8cm, clip, width=6cm]{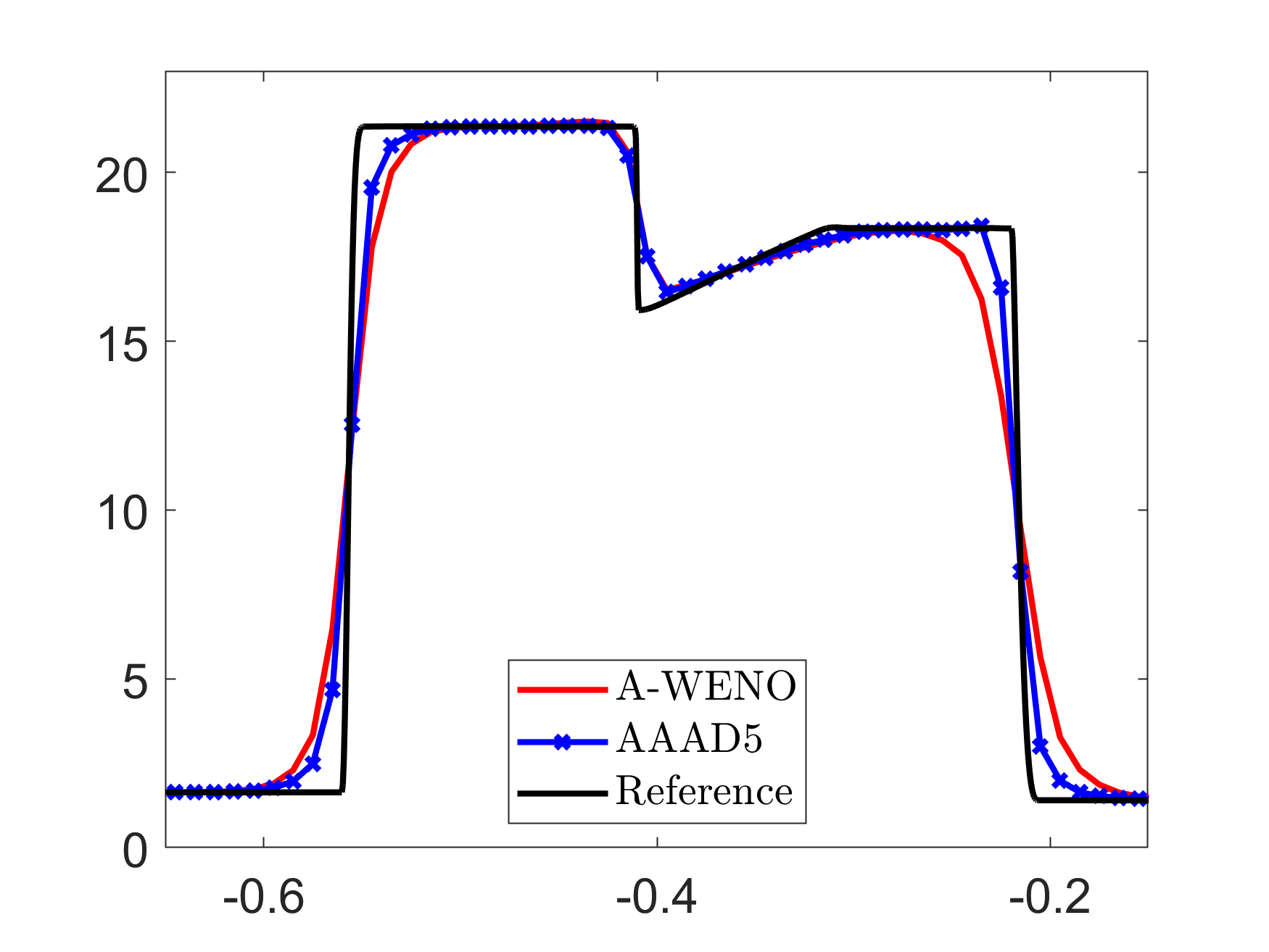}}
\caption{\sf Example 4: Density $\rho$ computed by the A-WENO and AAAD5 schemes (left) and zoom at $x\in[-0.65,-0.15]$ (right).
\label{fig3b}}
\end{figure}

We compute the numerical solution until the final time $t=3$ by the studied AAAD2 (with the adaptation constant $\texttt C=0.15$) and
AAAD5 (with the adaptation constant $\texttt C=0.05$) schemes on a uniform mesh with $\dx=1/100$. The obtained results are presented in
Figures \ref{fig3a}--\ref{fig3b} along with a reference solution computed by the CU scheme on a much finer mesh with $\dx=1/2000$. One can
see that both AAAD2 and AAAD5 achieve improved resolution of the contact discontinuities.

\begin{figure}[ht!]
\centerline{\includegraphics[trim=0.8cm 0.2cm 0.9cm 0.8cm, clip, width=5.cm]{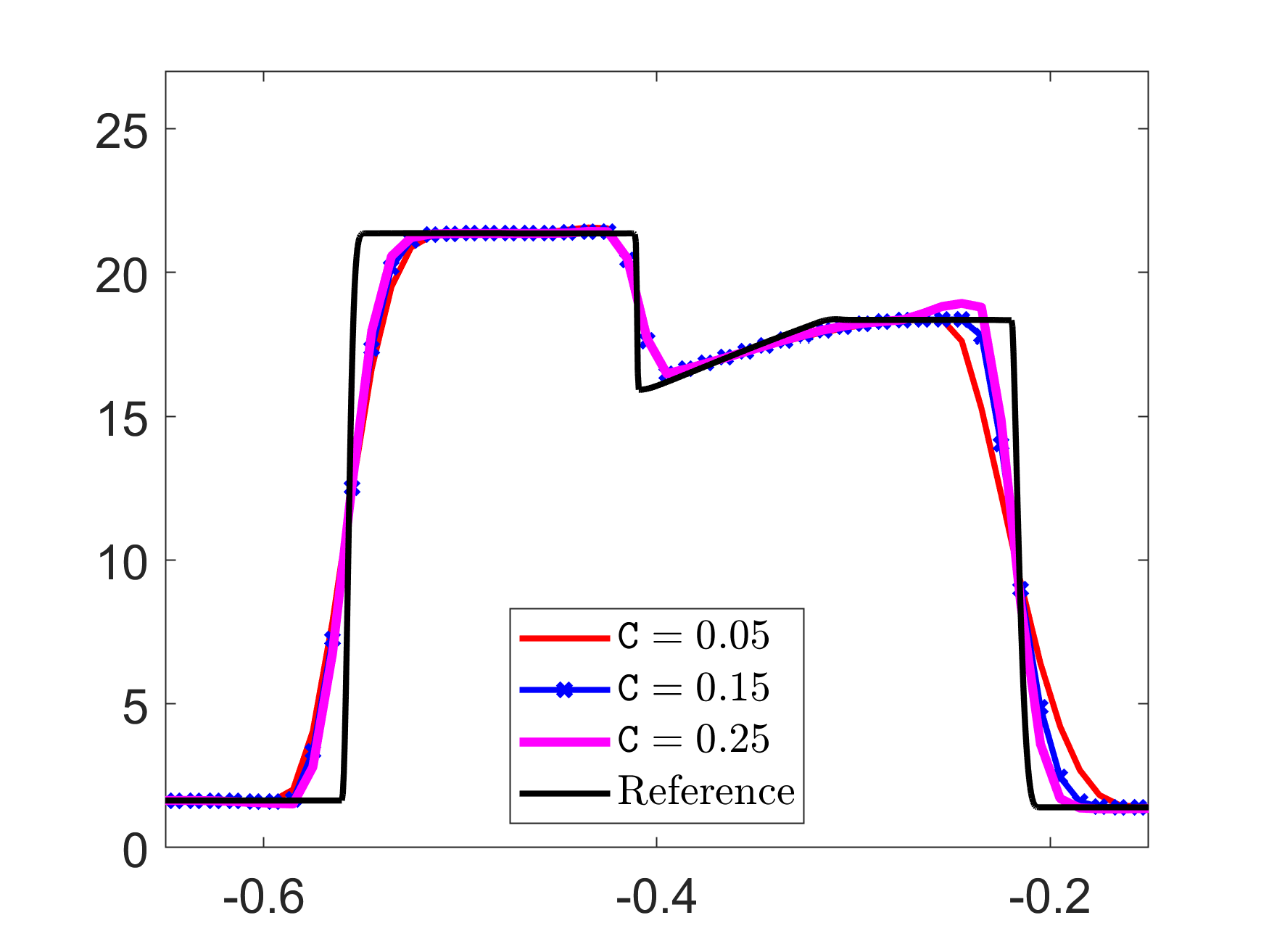}\hspace{0.3cm}
            \includegraphics[trim=0.8cm 0.2cm 0.9cm 0.8cm, clip, width=5.cm]{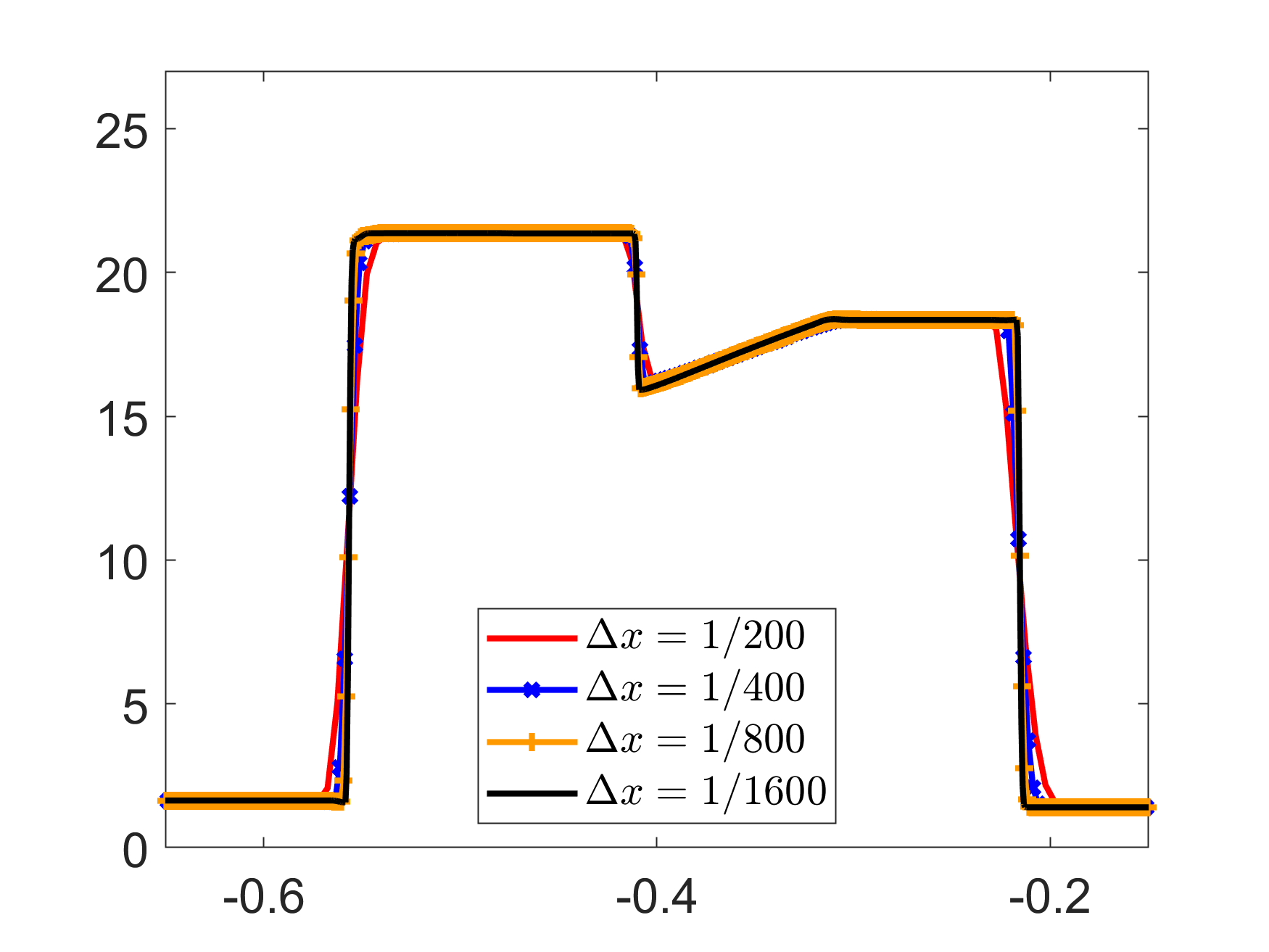}\hspace{0.3cm}
            \includegraphics[trim=0.8cm 0.2cm 0.9cm 0.8cm, clip, width=5.cm]{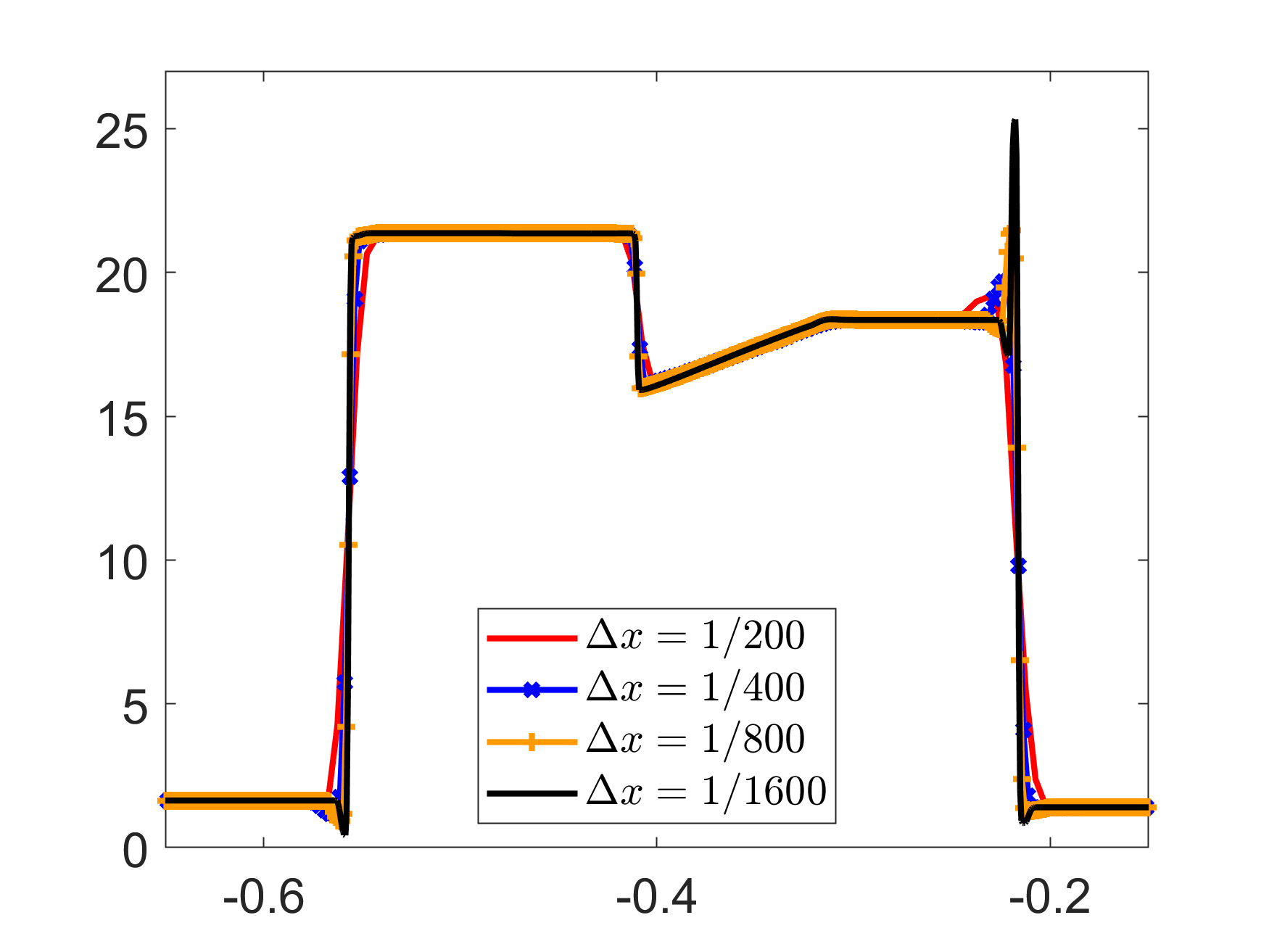}}
\caption{\sf Example 4: Density $\rho$ (zoom at $x\in[-0.65,-0.15]$) computed by the AAAD2 scheme with $\texttt C=0.05$, $0.15$, and $0.25$
on the uniform mesh with $\dx=1/100$ (left) and with $\texttt C=0.15$ (middle) and $\texttt C=0.25$ (right) on finer uniform meshes with
$\dx=1/200$, $1/400$, $1/800$, and $1/1600$.\label{fig4c}}
\end{figure}

We use this example to demonstrate the adaption coefficient tuning procedure. To this end, we first adjust $\texttt C$ on a coarse mesh (so
that the tuning process is computationally inexpensive) and then use the selected value for high-resolution computations on finer meshes. We
compute the numerical solutions by the AAAD2 scheme on a coarse mesh with $\dx=1/100$ using $\texttt C=0.05$, $0.15$, and $0.25$, and
present the results in Figure \ref{fig4c} (left), where one can see that $\texttt C=0.15$ seems to be a reasonable choice as the larger
value $\texttt C=0.25$ produces a small oscillation. We then take $\texttt C=0.15$, compute the corresponding numerical solutions on finer
meshes with $\dx=1/200$, $1/400$, $1/800$, and $1/1600$, and plot the obtained results in Figure \ref{fig4c} (middle). As one can see, the
computed solutions remain stable and well resolved as the mesh is refined. We also repeat the same finer mesh computations, but using the
larger value $\texttt C=0.25$. We plot the obtained results in Figure \ref{fig4c} (right), where one can see that the magnitude of
oscillations increases as $\dx$ shrinks, which clearly demonstrates that the coarse-mesh tuning process works in a robust way.

\subsubsection*{Example 5---Lax Problem}
In this example taken from \cite{Lax73}, we consider the Riemann problem with the initial conditions,
\begin{equation*}
(\rho,u,p)(x,0)=\begin{cases}(0.445,0.31061,8.928),&x<0,\\(0.5,0,0.571),&x\ge0,\end{cases}
\end{equation*}
prescribed in the interval $[-5,5]$ subject to the free boundary conditions.

\begin{figure}[ht!]
\centerline{\includegraphics[trim=1.0cm 0.3cm 0.9cm 0.6cm, clip, width=6.cm]{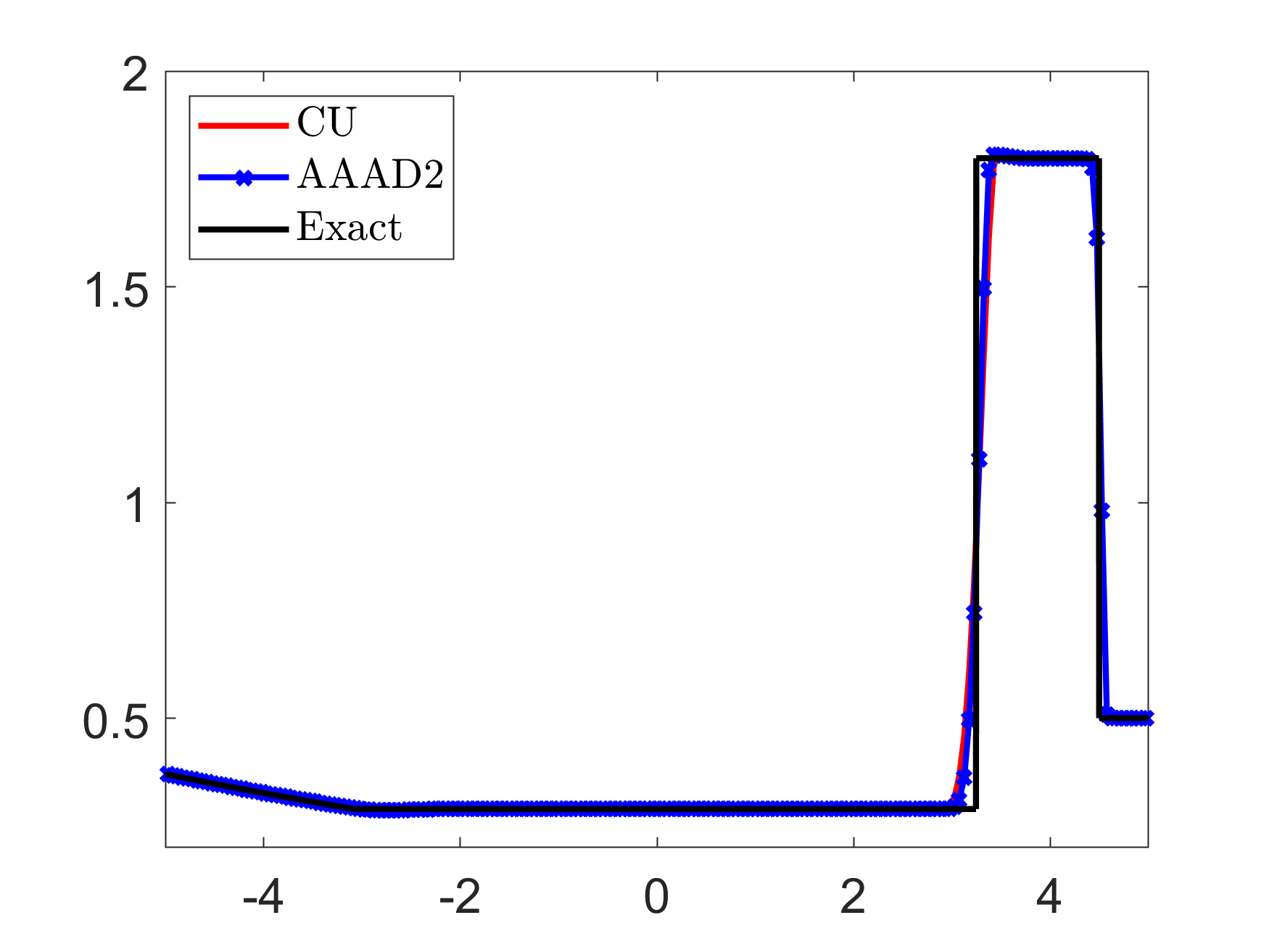}\hspace{1cm}
            \includegraphics[trim=1.0cm 0.3cm 0.9cm 0.6cm, clip, width=6.cm]{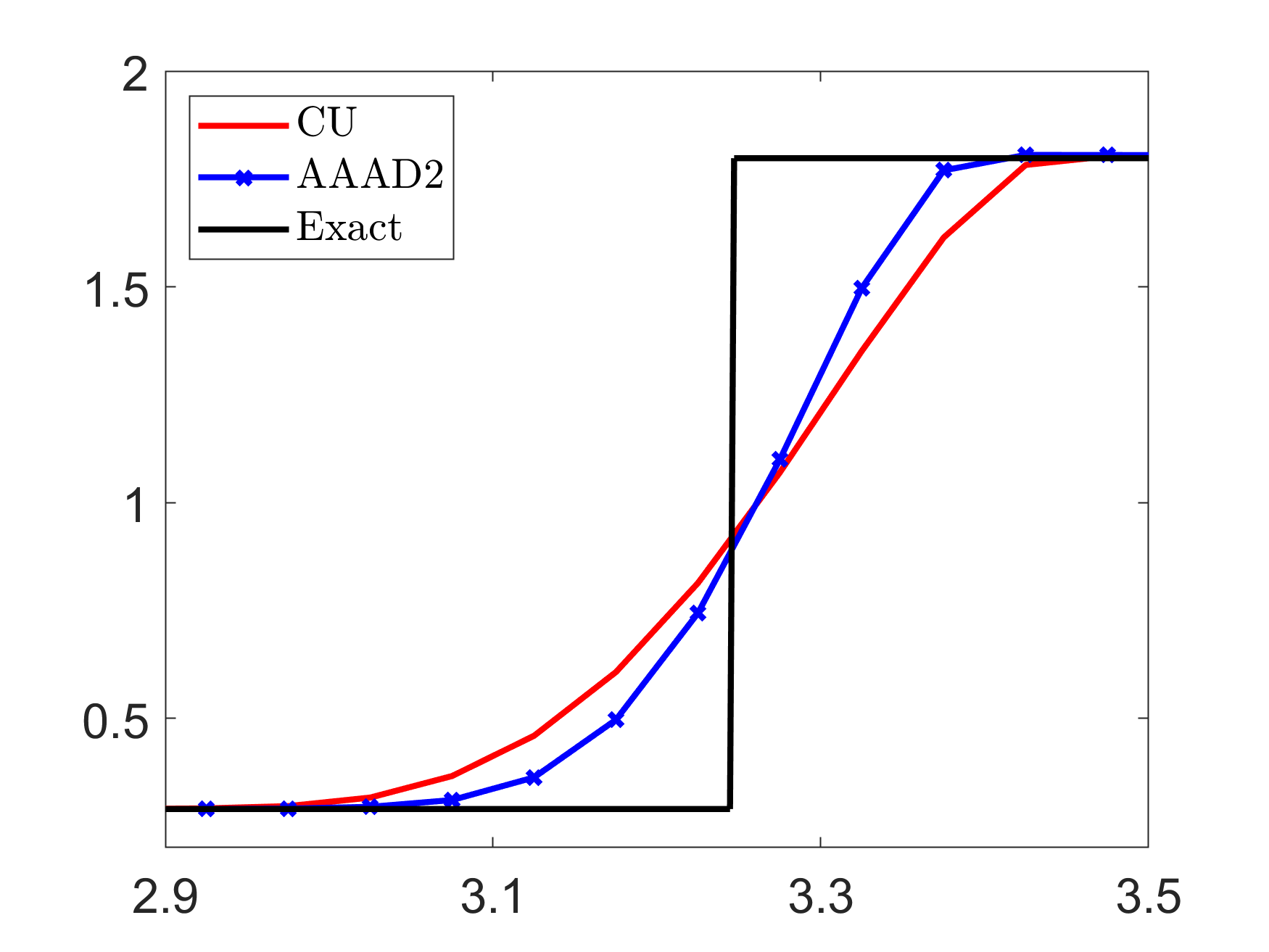}}
\caption{\sf Example 5: Density $\rho$ computed by the CU and AAAD2 schemes (left) and zoom at $x\in[2.9,3.5]$ (right).\label{fig2a}}
\end{figure}

\begin{figure}[ht!]
\centerline{\includegraphics[trim=1.0cm 0.3cm 0.9cm 0.6cm, clip, width=6.cm]{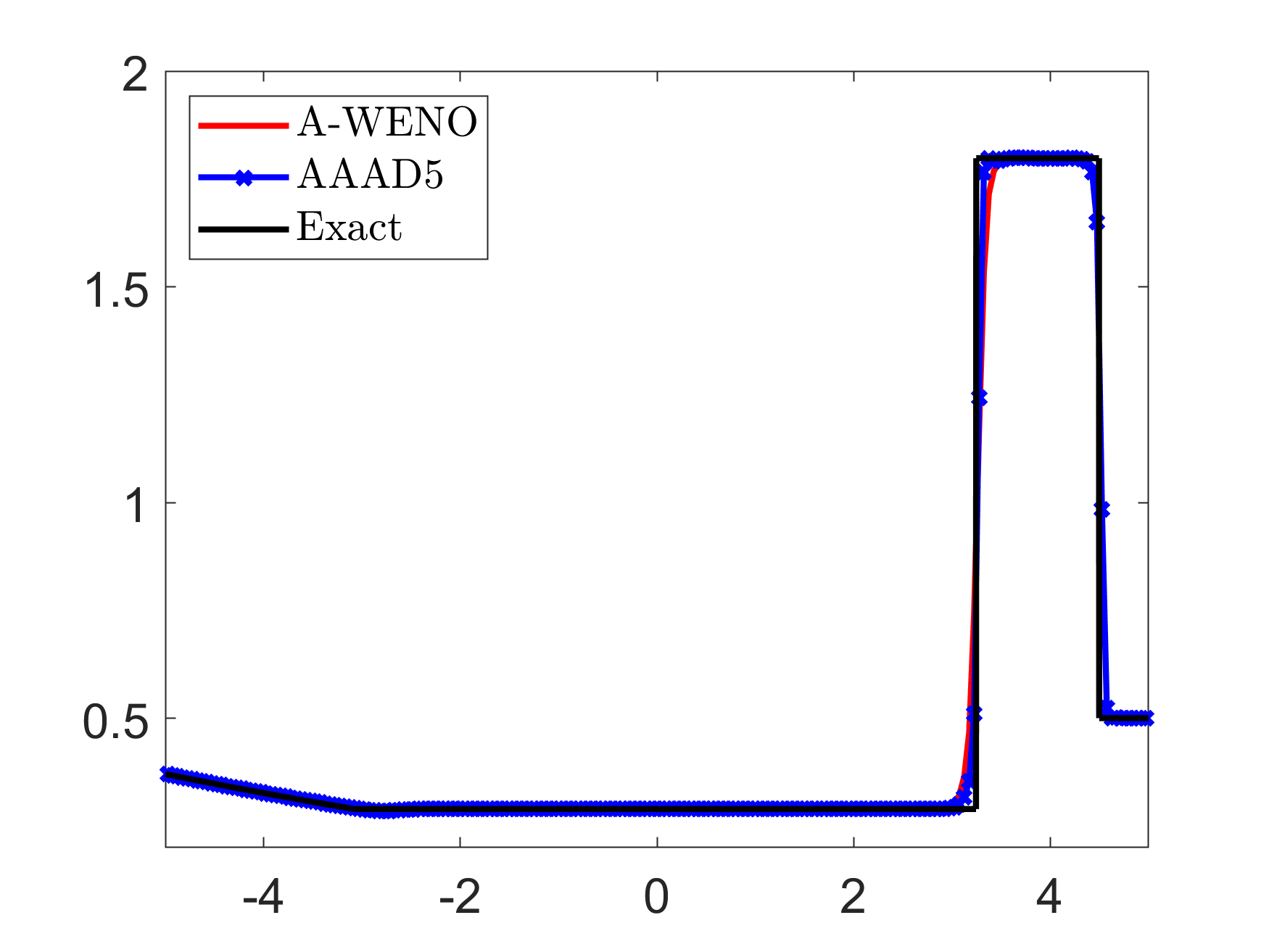}\hspace{1cm}
            \includegraphics[trim=1.0cm 0.3cm 0.9cm 0.6cm, clip, width=6.cm]{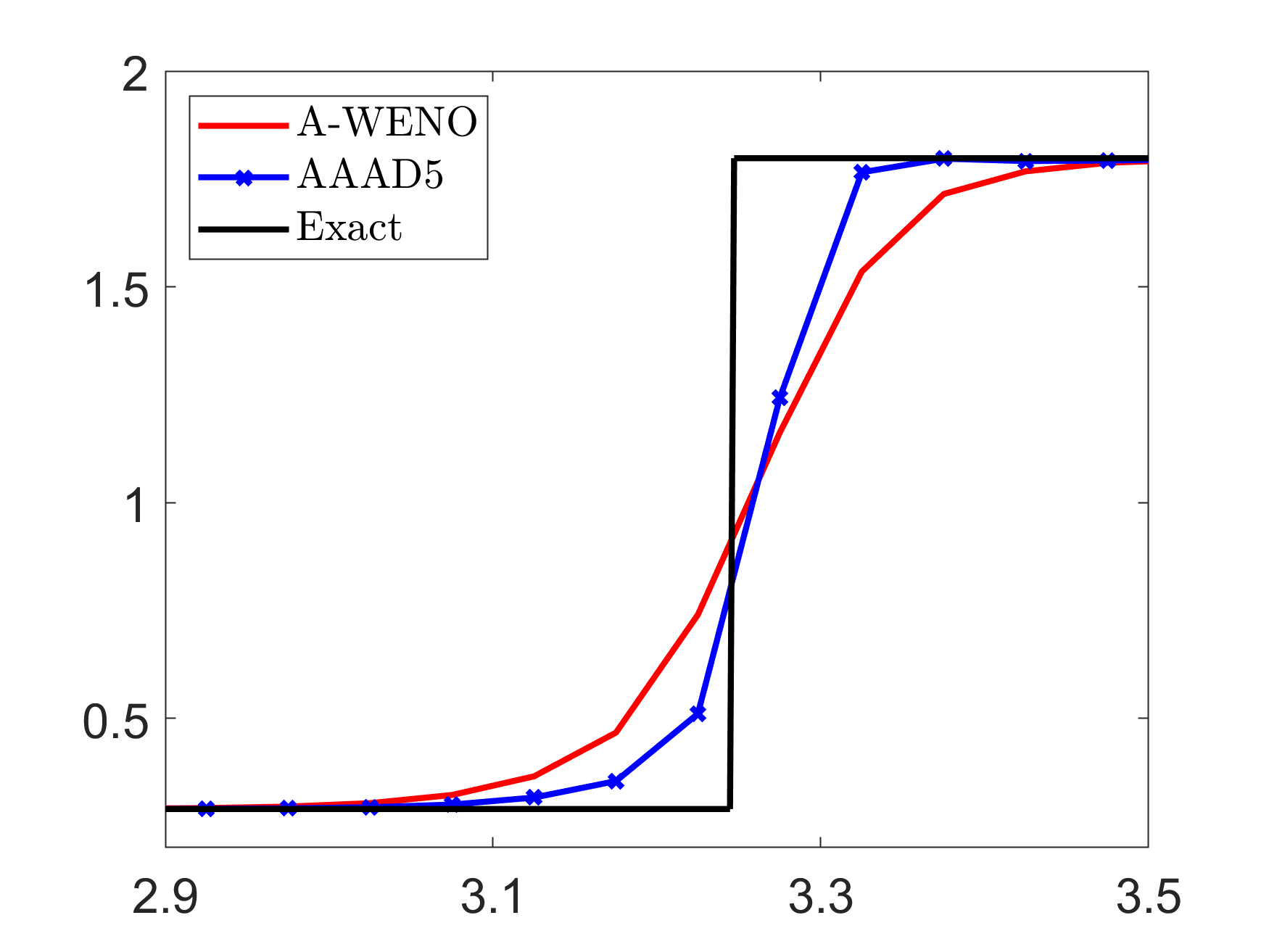}}
\caption{\sf Example 5: Density $\rho$ computed by the A-WENO and AAAD5 schemes (left) and zoom at $x\in[2.9,3.5]$ (right).\label{fig2b}}
\end{figure}

We compute the numerical solutions until the final time $t=1.3$ by the studied AAAD2 (with the adaptation constant $\texttt C=0.1$) and
AAAD5 (with the adaptation constant $\texttt C=0.5$) schemes on a uniform mesh with $\dx=1/20$. We plot the results in Figures
\ref{fig2a}--\ref{fig2b} together with the exact solution. As can be seen, both adaptive schemes provide improved resolution of the contact
wave compared with their non-adaptive counterparts.

\subsubsection*{Example 6---Blast Wave Problem}
In the last 1-D example, we consider a strong-shock interaction problem from \cite{Woodward88}, which is considered in the interval $[0,1]$
with the solid wall boundary conditions and subject to the following initial conditions:
\begin{equation*}
(\rho,u,p)\Big|_{(x,0)}=\begin{cases}(1,0,1000),&x<0.1,\\(1,0,0.01),&0.1\le x\le0.9,\\(1,0,100),&x>0.9.\end{cases}
\end{equation*}

\begin{figure}[ht!]
\centerline{\includegraphics[trim=1.0cm 0.3cm 0.8cm 0.8cm, clip, width=6.cm]{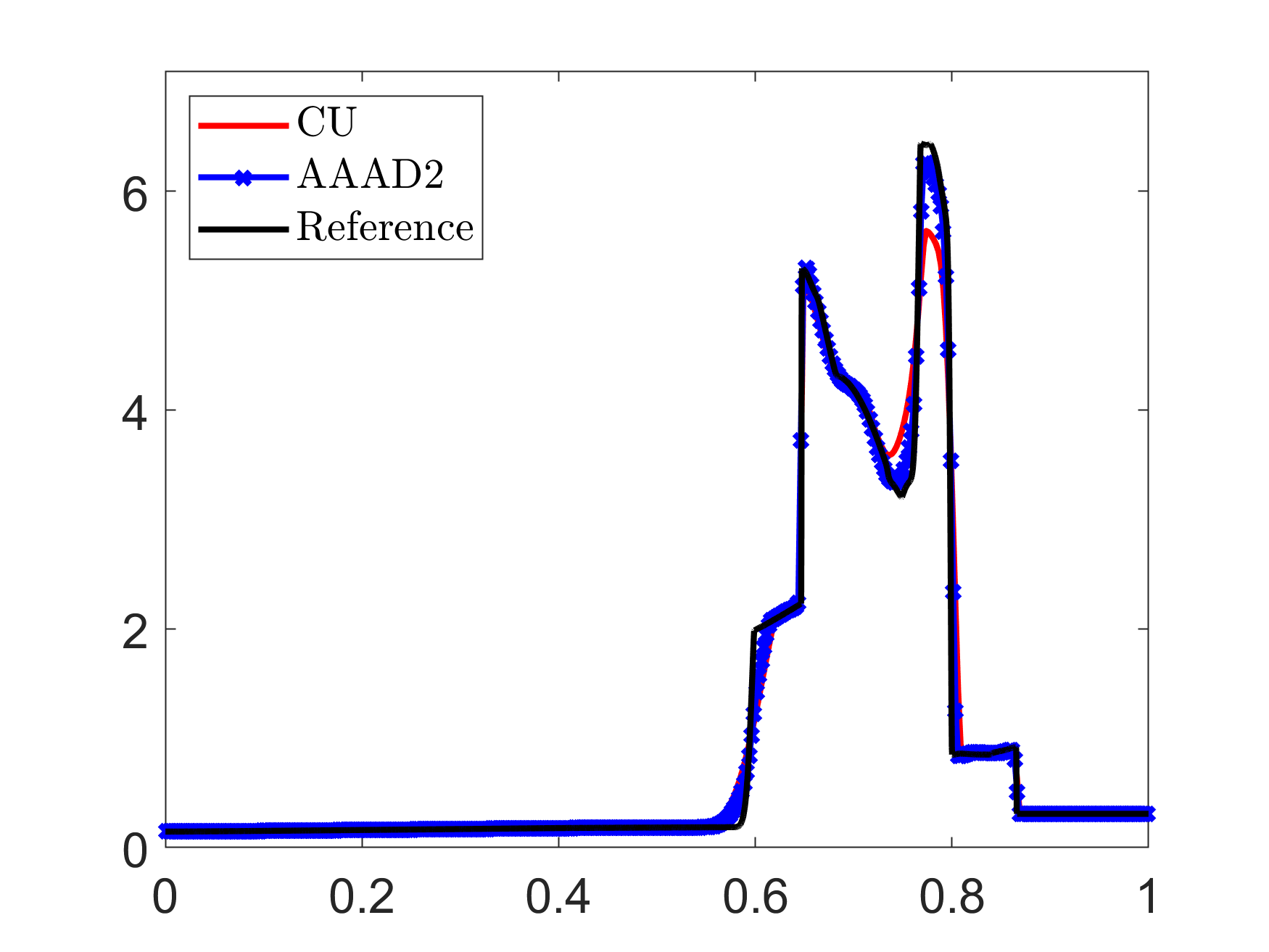}\hspace{1cm}
            \includegraphics[trim=1.0cm 0.3cm 0.8cm 0.8cm, clip, width=6.cm]{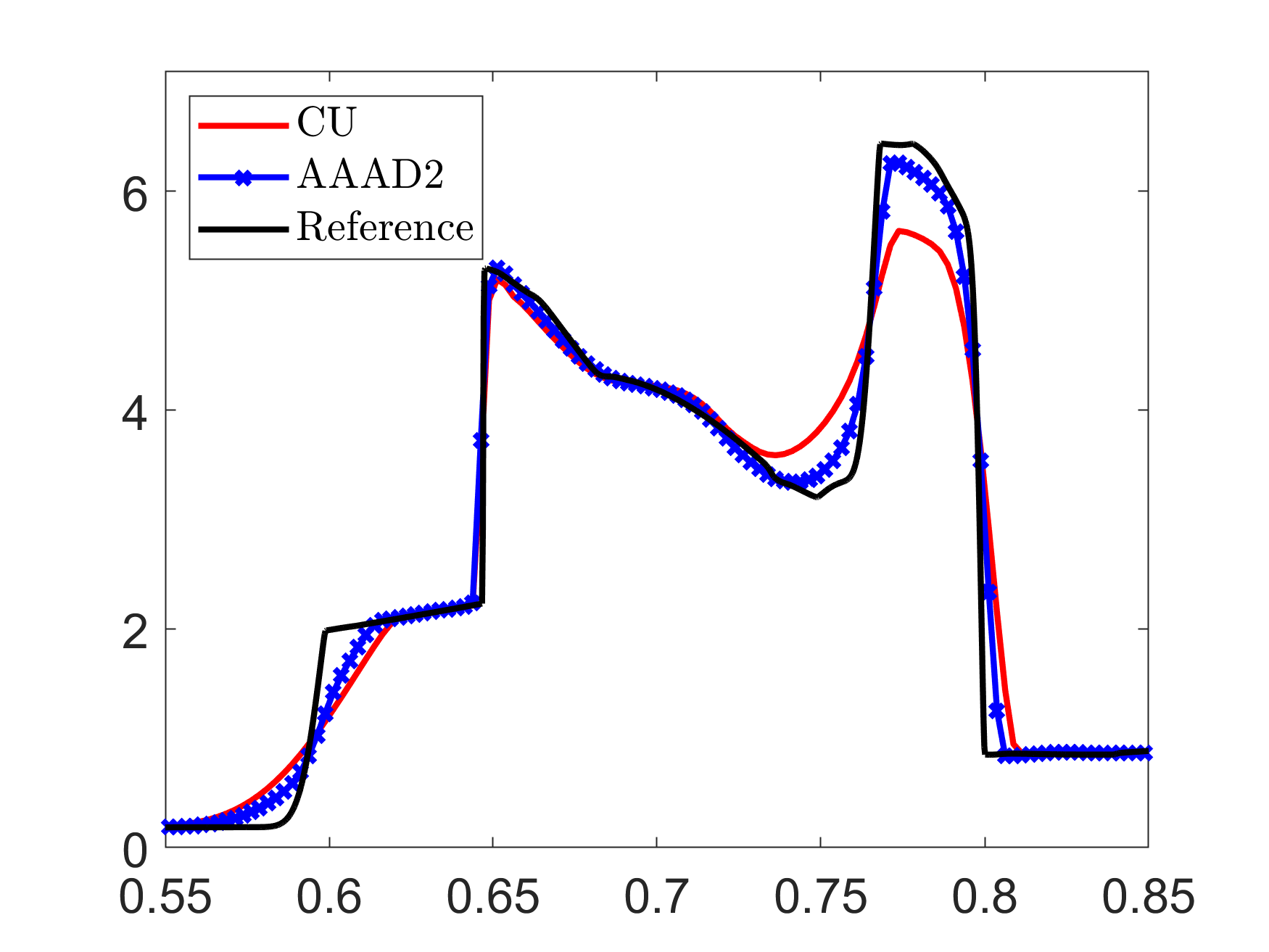}}
\caption{\sf Example 6: Density $\rho$ computed by the CU and AAAD2 schemes with $\dx=1/400$ (left) and zoom at $x\in[0.55,0.85]$ (right).
\label{fig3aa}}
\end{figure}

\begin{figure}[ht!]
\centerline{\includegraphics[trim=1.0cm 0.3cm 0.8cm 0.8cm, clip, width=6.cm]{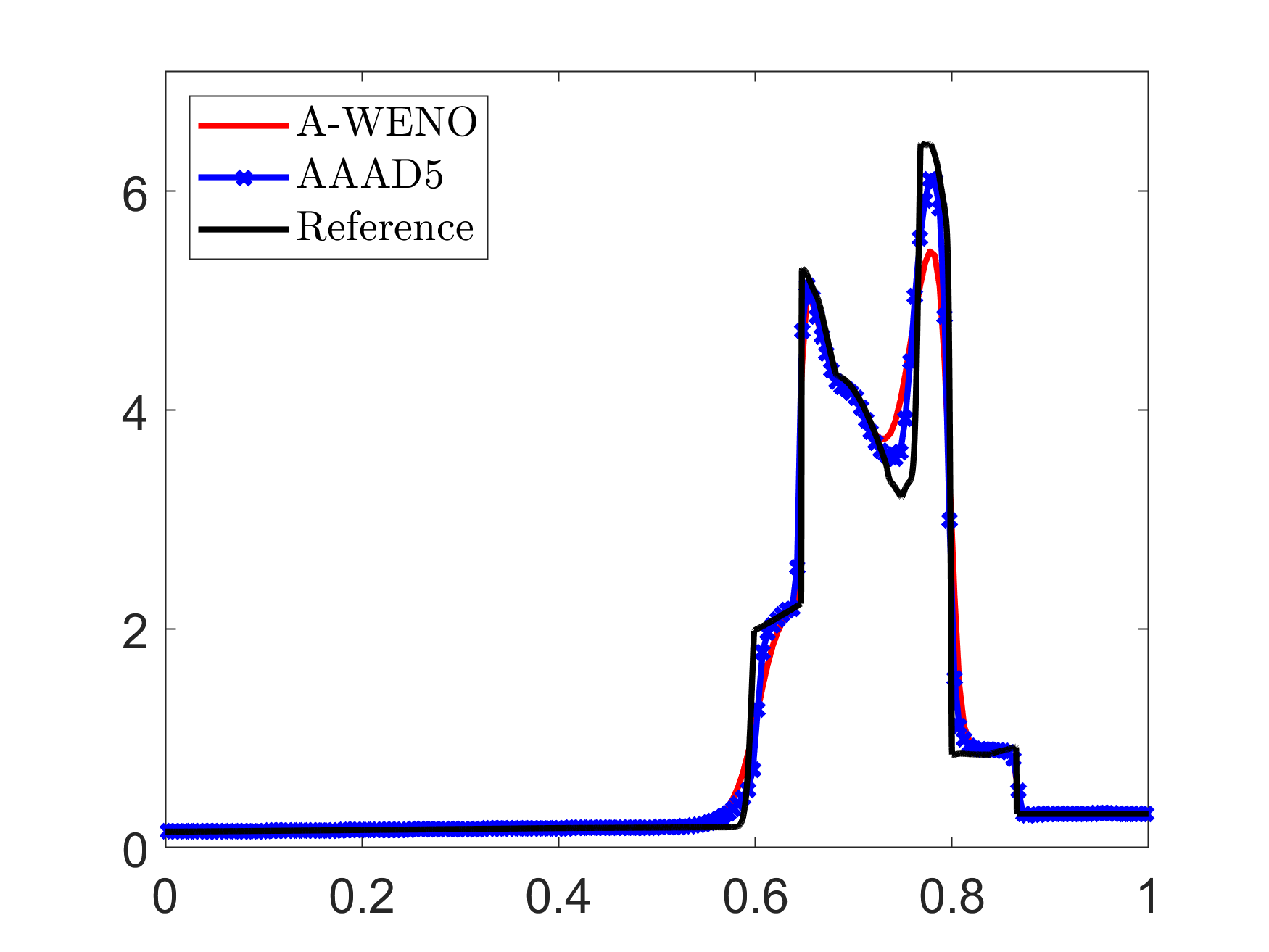}\hspace{1cm}
            \includegraphics[trim=1.0cm 0.3cm 0.8cm 0.8cm, clip, width=6.cm]{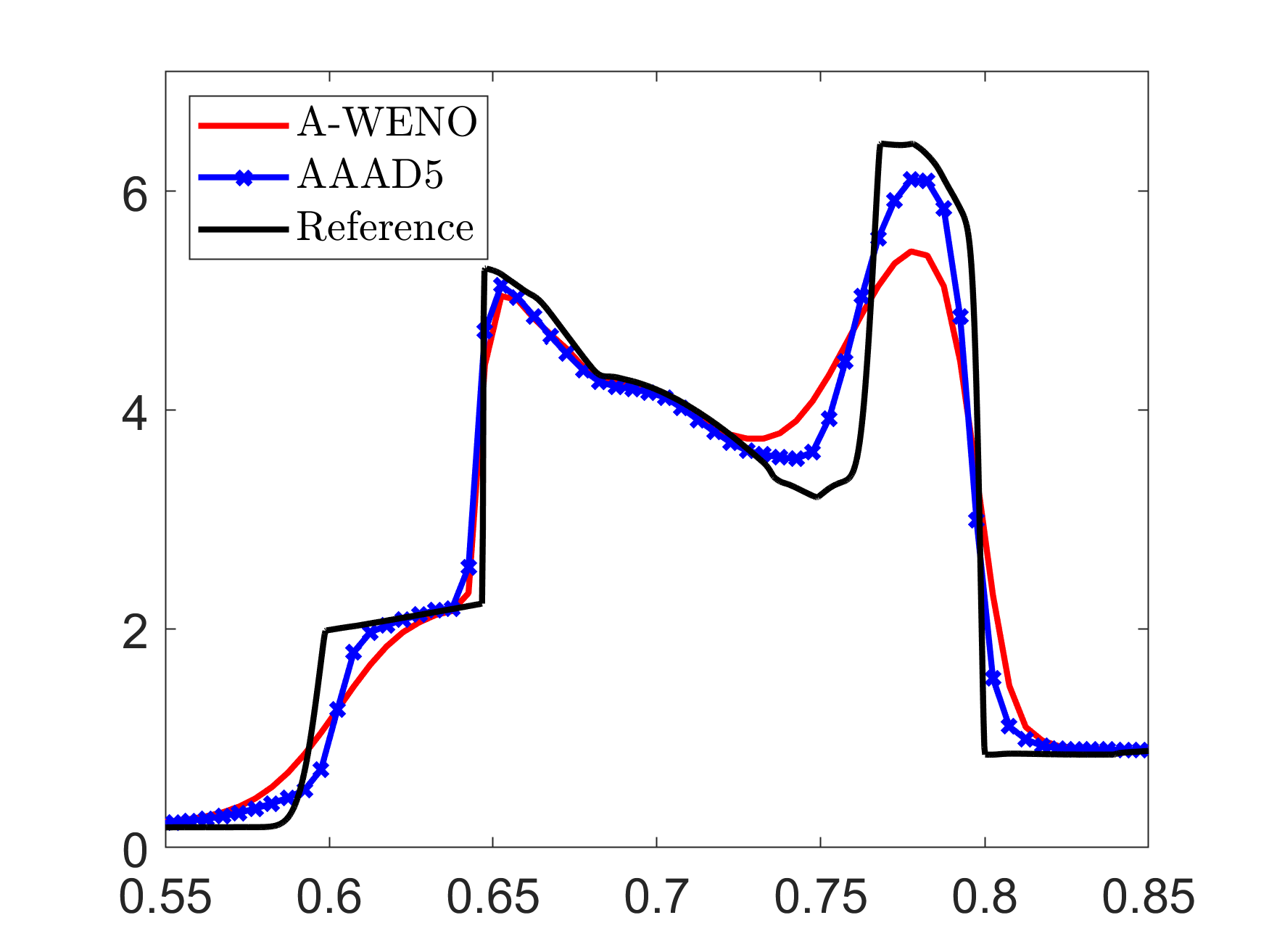}}
\caption{\sf Example 6: Density $\rho$ computed by the A-WENO and AAAD5 schemes with $\dx=1/200$ (left) and zoom at $x\in[0.55,0.85]$
(right).\label{fig3bb}}
\end{figure}

We compute the numerical solutions until the final time $t=0.038$ by the studied AAAD2 scheme (with the adaptation constant
$\texttt C=0.55$) on a uniform mesh with $\dx=1/400$ and AAAD5 scheme (with the adaptation constant $\texttt C=0.5$) on a coarser uniform
mesh with $\dx=1/200$. The obtained results are presented in Figures \ref{fig3aa}--\ref{fig3bb} along with a reference solution computed by
the CU scheme on a much finer mesh with $\dx=1/4000$. One can see that both AAAD2 and AAAD5 schemes achieve substantial improvement in the
quality of the computed solutions.

\subsection{Two-Dimensional Examples}\label{sec42}
We now turn to the 2-D Euler equations of gas dynamics. In Examples 7--13, we take $\gamma=1.4$, while in Example 14, we take $\gamma=5/3$.

\subsubsection*{Example 7---2-D Accuracy Test}
In the first 2-D example taken from \cite{BD2013,Shu1998}, we consider the following smooth initial data:
\begin{equation*}
\begin{aligned}
&\rho(x,y,0)=\bigg(1-\frac{(\gamma-1)\kappa^2}{2\gamma}\bigg)^{\frac{1}{\gamma-1}},\quad p(x,y,0)=\rho^\gamma(x,y,0),\\
&u(x,y,0)=1-\kappa y, \quad v(x,y,0)=1+\kappa x, \quad \kappa = \frac{5}{2 \pi} e^{\frac{1-x^2-y^2}{2}} 
\end{aligned}
\end{equation*}
subject to the periodic boundary conditions in the computational domain $[-10,10]\times[-10,10]$. The exact solution of this initial value
problem is given by $\mU(x,y,t)=\mU(x-t,y-t,0)$.

We compute the numerical solution until the final time $t=0.1$ using the studied AAAD2 and AAAD5 schemes (both with the adaptation constant
$\texttt C=0.1$) on a sequence of uniform meshes with $\dx=\dy=1/10$, $1/20$, and $1/40$, and then measure the $L^1$-errors and the
corresponding experimental convergence rates for the density. The obtained results are presented in Table \ref{tab52}, where one can see
that the expected orders of accuracy have been achieved by both of the studied schemes.
\begin{table}[ht!]
\centering
\begin{tabular}{|c|cc|cc|}
\hline
\multirow{2}{*}{$\dx=\dy$}&\multicolumn{2}{c|}{AAAD2}&\multicolumn{2}{c|}{AAAD5}\\
\cline{2-5}
&Error&Rate&Error&Rate\\
\hline
$1/10$&3.62e-03&--- &3.73e-05&--- \\
$1/20$&8.18e-04&2.15&1.14e-06&5.04\\
$1/40$&1.69e-04&2.27&2.44e-08&5.55\\
\hline
\end{tabular}
\caption{\sf Example 7: The $L^1$-errors of the density $\rho$ and experimental convergence rates for the AAAD2 and AAAD5 schemes.
\label{tab52}}
\end{table}

\begin{remark}
As in Example 1, we use smaller time steps with $\dt\sim\min\{(\dx)^{5/3},(\dy)^{5/3}\}$ in order to achieve fifth order of accuracy by the
AAAD5 scheme.
\end{remark}

\subsubsection*{Example 8---Explosion Problem}
In this example taken from \cite{Liska03}, we consider the explosion problem. The initial conditions,
\begin{equation*}
(\rho(x,y,0),u(x,y,0),v(x,y,0),p(x,y,0))=\begin{cases}(1,0,0,1),&x^2+y^2<0.16,\\(0.125,0,0,0.1),&\mbox{otherwise},\end{cases}
\end{equation*}
are prescribed in $[-1.5,1.5]\times[-1.5,1.5]$ subject to the free boundary conditions. The solution develops circular shock and contact
waves, with the latter one being unstable. 

\begin{figure}[ht!]
\centerline{\includegraphics[trim=0cm 0cm 0cm 0cm, clip, width=\linewidth]{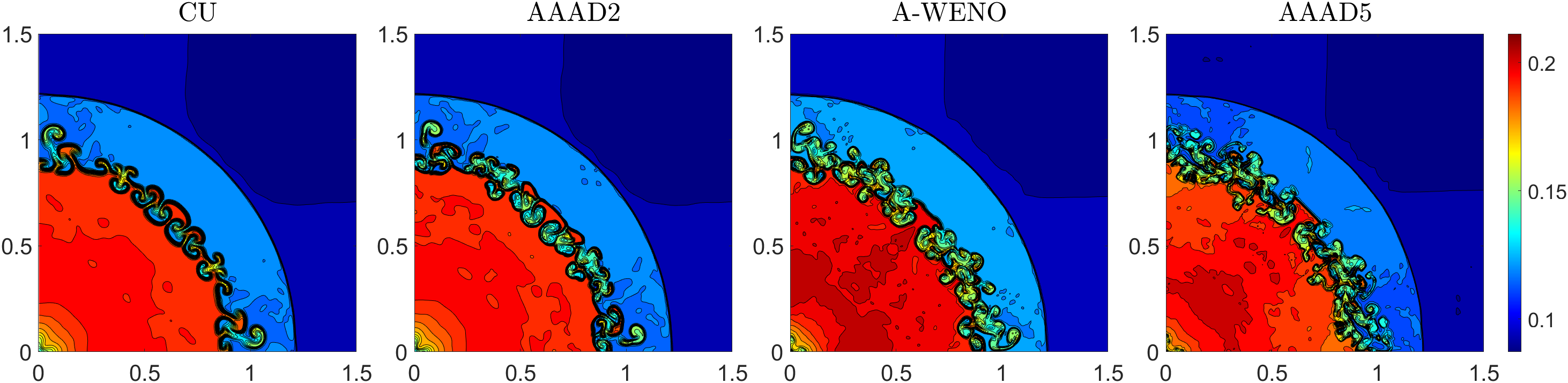}}
\caption{\sf Example 8: Density $\rho$ computed by the CU, AAAD2, A-WENO, and AAAD5 schemes.\label{fig8}}
\end{figure}

We compute the numerical solutions until the final time $t=3.2$ by the studied AAAD2 (with the adaptation constant $\texttt C=0.03$) and
AAAD5 (with the adaptation constant $\texttt C=0.02$) schemes on a uniform mesh with $\dx=\dy=3/800$. The obtained results are presented in
Figure \ref{fig8}, where one can clearly observe noticeable differences between the CU and A-WENO schemes and their AAAD counterparts: the
adaptive solutions exhibit substantially more small-scale structures and a wider mixing layer, while the main shock remains stable.

\subsubsection*{Example 9---2-D Riemann Problem (Configuration 3)}
In this example, we consider Configuration 3 of the 2-D Riemann problems from \cite{Kurganov02} with the initial conditions
\begin{equation*}
(\rho,u,v,p)(x,y,0)=\begin{cases}(1.5,0,0,1.5),&x>1,~y>1,\\(0.5323,1.206,0,0.3),&x<1,~y>1,\\(0.138,1.206,1.206,0.029),&x<1,~y<1,\\
(0.5323,0,1.206,0.3),&x>1,~y<1,\end{cases}
\end{equation*}
prescribed in $[0,1.2]\times[0,1.2]$ subject to the free boundary conditions.

\begin{figure}[ht!]
\centerline{\includegraphics[trim=0cm 0cm 0cm 0cm, clip, width=\linewidth]{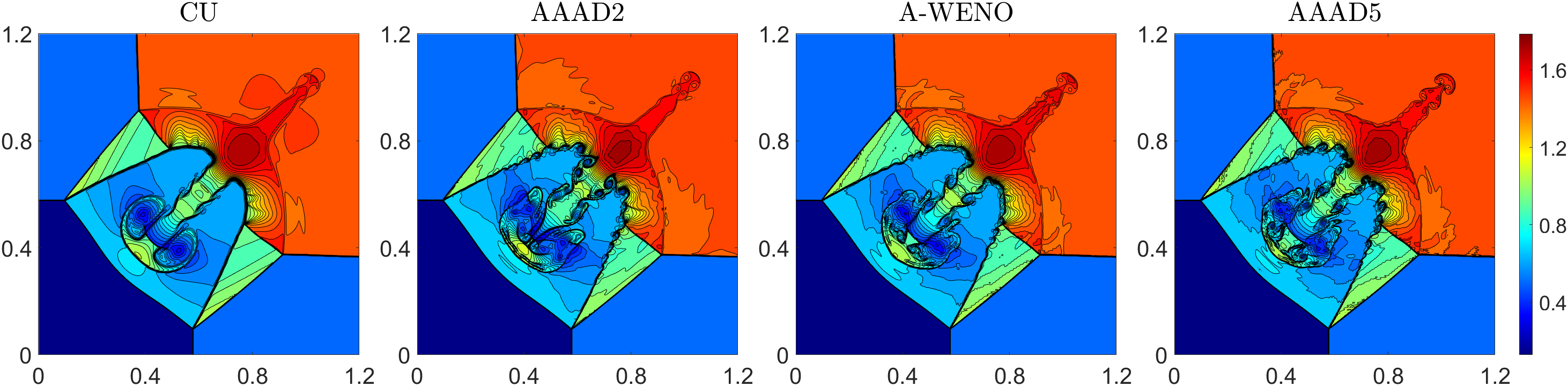}}
\caption{\sf Example 9: Density $\rho$ computed by the CU, AAAD2, A-WENO, and AAAD5 schemes.\label{fig9}}
\end{figure}

We compute the numerical solution until the final time $t=1$ by the studied AAAD2 (with the adaptation constant $\texttt C=0.04$) and AAAD5
(with the adaptation constant $\texttt C=0.02$) schemes on a uniform mesh with $\dx=\dy=1/500$ and plot the results in Figure \ref{fig9}.
As one can see, the AAAD schemes outperform their non-adaptive counterparts in capturing the sideband instability of the jet, in particular,
in the zones of strong along-jet velocity shear and near the jet neck.

\begin{figure}[ht!]
\centerline{\includegraphics[trim=0cm 0cm 0cm 0cm, clip, width=\linewidth]{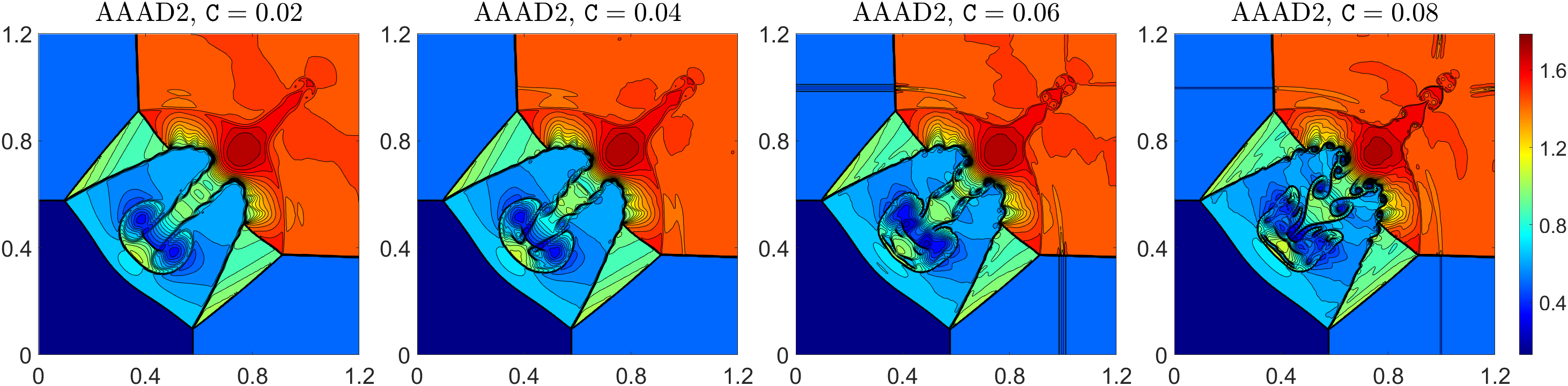}}
\caption{\sf Example 9: Density $\rho$ computed by the AAAD2 scheme with $\dx=\dy=3/1000$ for $\texttt C=0.02$, $0.04$, $0.06$, and $0.08$.
\label{fig9a}}
\end{figure}

\begin{figure}[ht!]
\centerline{\includegraphics[trim=0cm 0cm 0cm 0cm, clip, width=\linewidth]{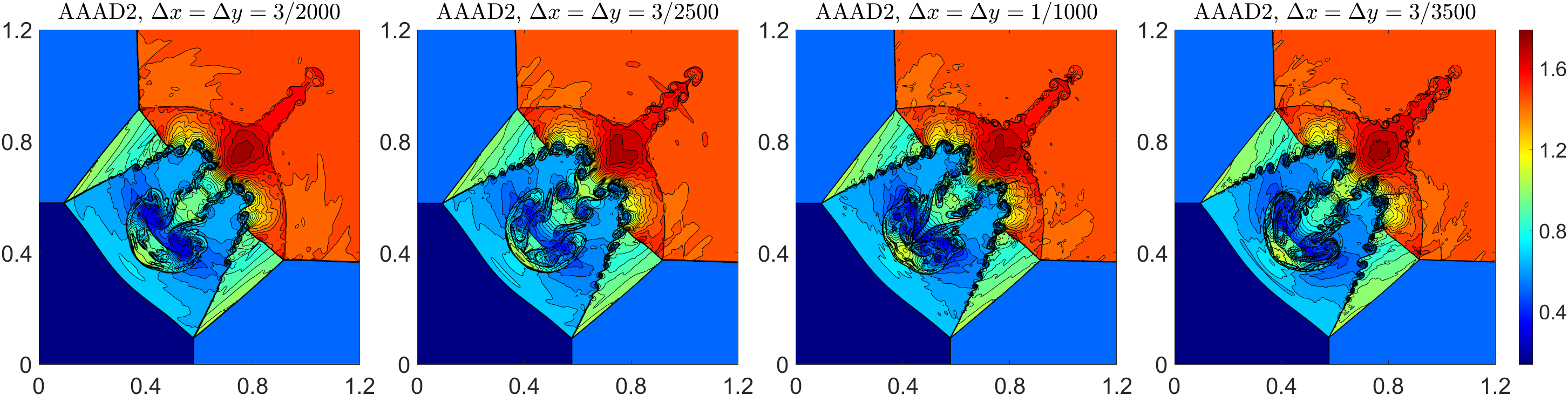}}
\caption{\sf Example 9: Density $\rho$ computed by the AAAD2 scheme on uniform meshes with $\dx=\dy=3/2000$, $3/2500$, $1/1000$, and
$3/3500$ with $\texttt C=0.04$.\label{fig9b}}
\end{figure}

We use this example to illustrate the tuning procedure for $\texttt C$ in the 2-D setting. To this end, we proceed as in the 1-D case: we
first tune $\texttt C$ on a coarse mesh and then use the selected values for fine-mesh computations. We begin by computing the coarse-mesh
numerical solutions with $\dx=\dy=3/1000$ using $\texttt C=0.02$, $0.04$, $0.06$, and $0.08$, and plot the results in Figure \ref{fig9a},
where one can see that $\texttt C=0.04$ seems to be a reasonable choice as the use of larger values of $\texttt C$ leads to the appearance
of certain unstable structures. We then verify that the value $\texttt C=0.04$ can be safely used on finer meshes with $\dx=\dy=3/2000$,
$3/2500$, $1/1000$, and $3/3500$; see Figure \ref{fig9b}. We have also run the AAAD2 scheme with $\texttt C=0.06$ on finer meshes, but when
the mesh is refined, the magnitude of oscillations grows and negative pressure values start appearing.

\subsubsection*{Example 10---2-D Riemann Problem (Configuration 6)}
In this example, we consider Configuration 6 of the 2-D Riemann problems taken from \cite{Kurganov02} with the initial conditions
\begin{equation*}
(\rho,u,v,p)(x,y,0)=\begin{cases}(1,0.75,-0.5,1),&x>0.5,~y>0.5,\\(2,0.75,0.5,1),&x<0.5,~y>0.5,\\(1,-0.75,0.5,1),&x<0.5,~y<0.5,\\
(3,-0.75,-0.5,1),&x>0.5,~y<0.5,\end{cases}
\end{equation*}
prescribed in $[0,1]\times[0,1]$ subject to the free boundary conditions.

\begin{figure}[ht!]
\centerline{\includegraphics[trim=0cm 0cm 0cm 0cm, clip, width=\linewidth]{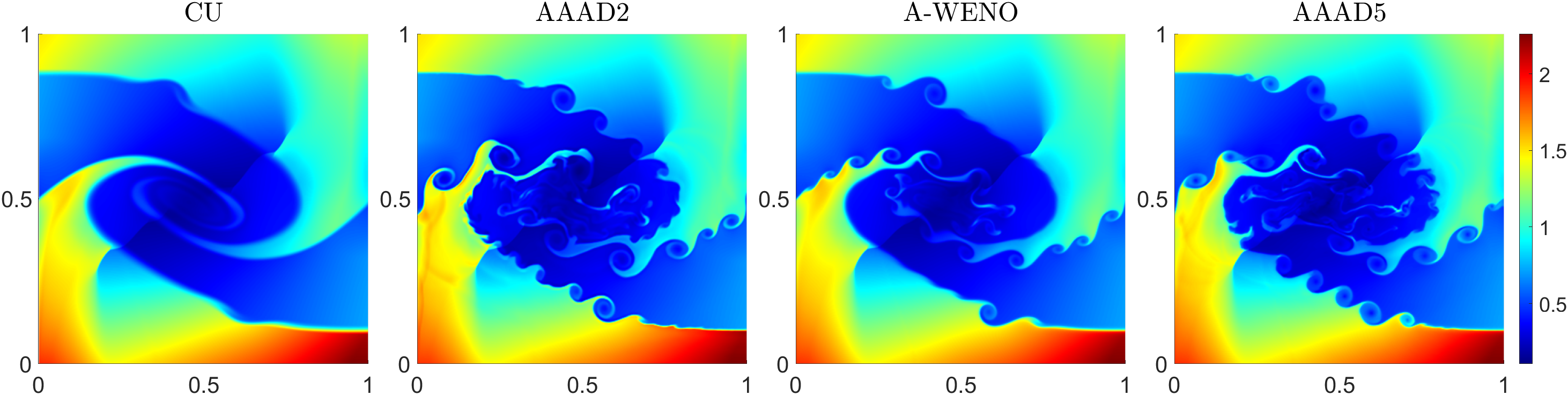}}
\caption{\sf Example 10: Density $\rho$ computed by the CU, AAAD2, A-WENO, and AAAD5 schemes.\label{fig10}}
\end{figure}

We compute the numerical solution until the final time $t=1$ by the studied AAAD2 (with the adaptation constant $\texttt C=0.05$) and AAAD5
(with the adaptation constant $\texttt C=0.02$) schemes on a uniform mesh with $\dx=\dy=1/400$ and plot the results in Figure \ref{fig10}.
As one can see, the AAAD schemes capture more intricate vortex structures than their non-adaptive counterparts, which indicates higher
resolution and reduced numerical dissipation of the proposed AAAD schemes. One can also observe that the second-order AAAD2 scheme achieves
higher resolution than the fifth-order A-WENO scheme. This demonstrates a powerful feature of the proposed AAAD strategy.

\subsubsection*{Example 11---2-D Riemann Problem (Configuration 12)}
In this example, we consider Configuration 12 of the 2-D Riemann problems taken from \cite{Kurganov02} with the initial conditions
\begin{equation*}
(\rho,u,v,p)\Big|_{(x,y,0)}=\begin{cases}(0.5313,0,0,0.4),&x>0.5,~y>0.5,\\(1,0.7276,0,1),&x<0.5,~y>0.5,\\(0.8,0,0,1),&x<0.5,~y<0.5,\\
(1,0,0.7276,1),&x>0.5,~y<0.5,\end{cases}
\end{equation*}
prescribed in the computational domain $[0,0.6]\times[0,0.6]$ subject to the free boundary conditions.

\begin{figure}[ht!]
\centerline{\includegraphics[trim=0cm 0cm 0cm 0cm, clip, width=\linewidth]{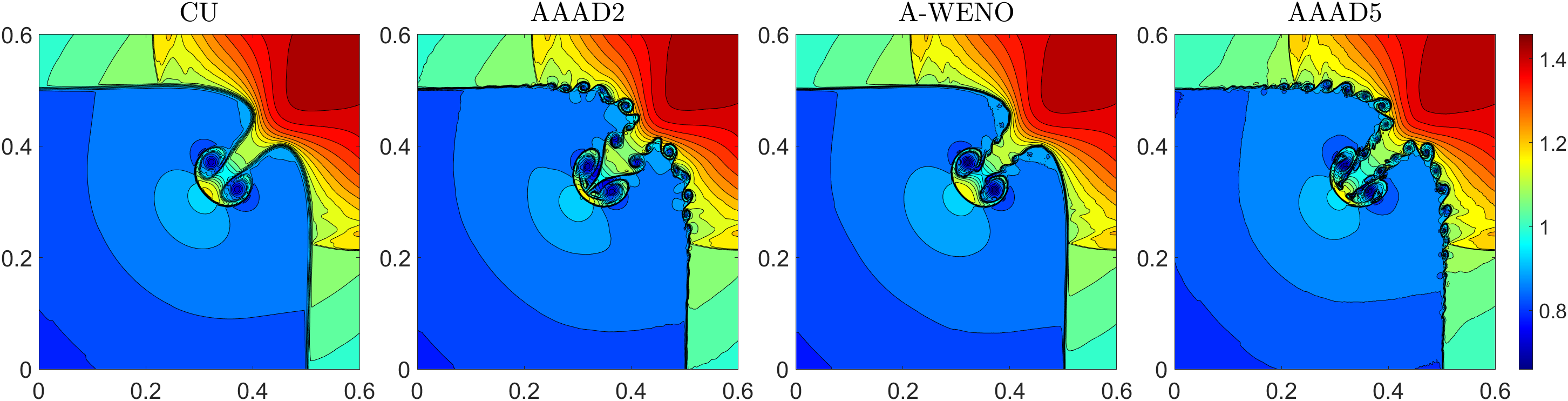}}
\caption{\sf Example 11: Density $\rho$ computed by the CU, AAAD2, A-WENO, and AAAD5 schemes.\label{fig10a}}
\end{figure}

We compute the numerical solution until the final time $t=1$ by the studied AAAD2 (with the adaptation constant $\texttt C=0.04$) and AAAD5
(with the adaptation constant $\texttt C=0.02$) schemes on a uniform mesh with $\dx=\dy=1/1000$ and plot the obtained results in Figure
\ref{fig10a}. One can see that the AAAD schemes resolve more vortices arising along the unstable contact surfaces than their non-adaptive
counterparts. As in the previous example, we notice that the second-order AAAD2 scheme clearly outperforms the fifth-order A-WENO scheme.

\subsubsection*{Example 12---Implosion Problem}
In this example, we consider the implosion problem taken from \cite{Liska03}; see also \cite{CCHKL_22,Kurganov07}. The initial conditions
\begin{equation*}
(\rho(x,y,0),u(x,y,0),v(x,y,0),p(x,y,0))=\begin{cases}(0.125,0,0,0.14),&|x|+|y|<0.15,\\(1,0,0,1),&\mbox{otherwise},\end{cases}
\end{equation*}
are prescribed in $[0,0.3]\times[0,0.3]$ subject to the solid boundary conditions. In this problem, a jet forms near the origin and
propagates along the diagonal $y=x$. Excessive numerical diffusion may smear the jet or alter its propagation speed, and therefore this
example is often used to assess numerical dissipation present in studied schemes.

\begin{figure}[ht!]
\centerline{\includegraphics[trim=0cm 0cm 0cm 0cm, clip, width=\linewidth]{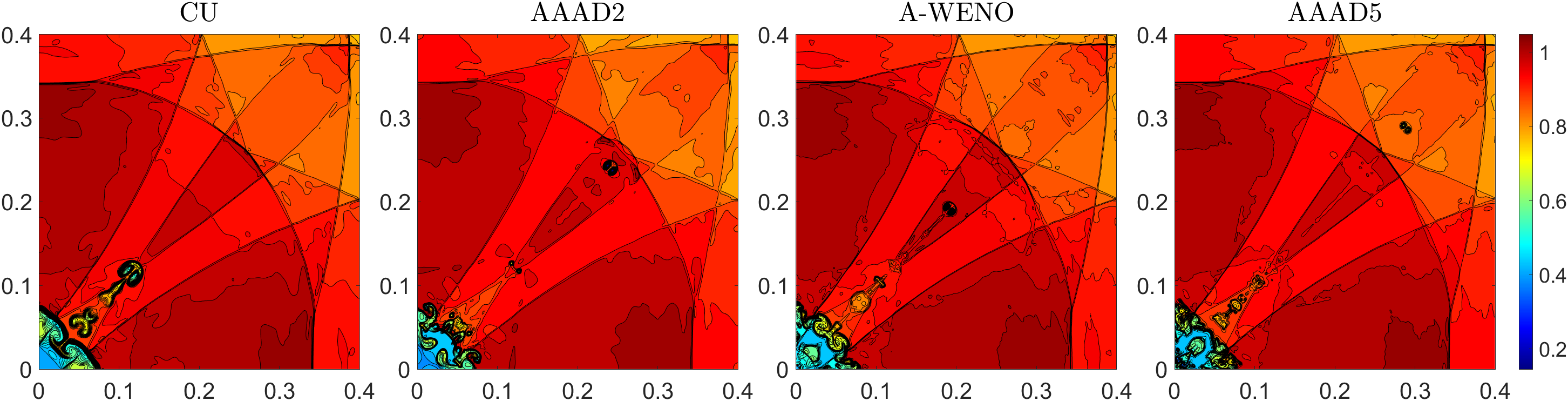}}
\caption{\sf Example 12: Density $\rho$ computed by the CU, AAAD2, A-WENO, and AAAD5 schemes.\label{fig11}}
\end{figure}

We compute the numerical solution until the final time $t=2.5$ by the studied AAAD2 (with the adaptation constant $\texttt C=0.05$) and
AAAD5 (with the adaptation constant $\texttt C=0.01$) schemes on a uniform mesh with $\dx=\dy=1/1500$. The obtained results are plotted in
Figure \ref{fig11}, where one can observe that the jets propagate substantially farther along the diagonal by the AAAD schemes: This
indicates that they are less diffusive than their non-adaptive counterparts. As in several previous examples, we notice that the
second-order AAAD2 scheme clearly outperforms the fifth-order A-WENO scheme.

\subsubsection*{Example 13---Kelvin-Helmholtz (KH) Instability}
In this example, we study the KH instability problem taken from \cite{Panuelos20}. We consider the initial data
\begin{equation}
\begin{aligned}
&(\rho(x,y,0),u(x,y,0))=\begin{cases}(1,-0.5+0.5e^{(y+0.25)/L}),&y<-0.25,\\(2,0.5-0.5e^{(-y-0.25)/L}),&-0.25<y<0,\\
(2,0.5-0.5e^{(y-0.25)/L}),&0<y<0.25,\\(1,-0.5+0.5e^{(0.25-y)/L}),&y>0.25,\end{cases}\\
&v(x,y,0)=0.01\sin(4\pi x),\quad p(x,y,0)\equiv1.5,
\end{aligned}
\label{4.1}
\end{equation}
in $[-0.5,0.5]\times[-0.5,0.5]$ subject to the periodic boundary conditions. In \eref{4.1}, $L=0.00625$ is a smoothing parameter.

\begin{figure}[ht!]
\centerline{\includegraphics[trim=0cm 0cm 0cm 0cm, clip, width=\linewidth]{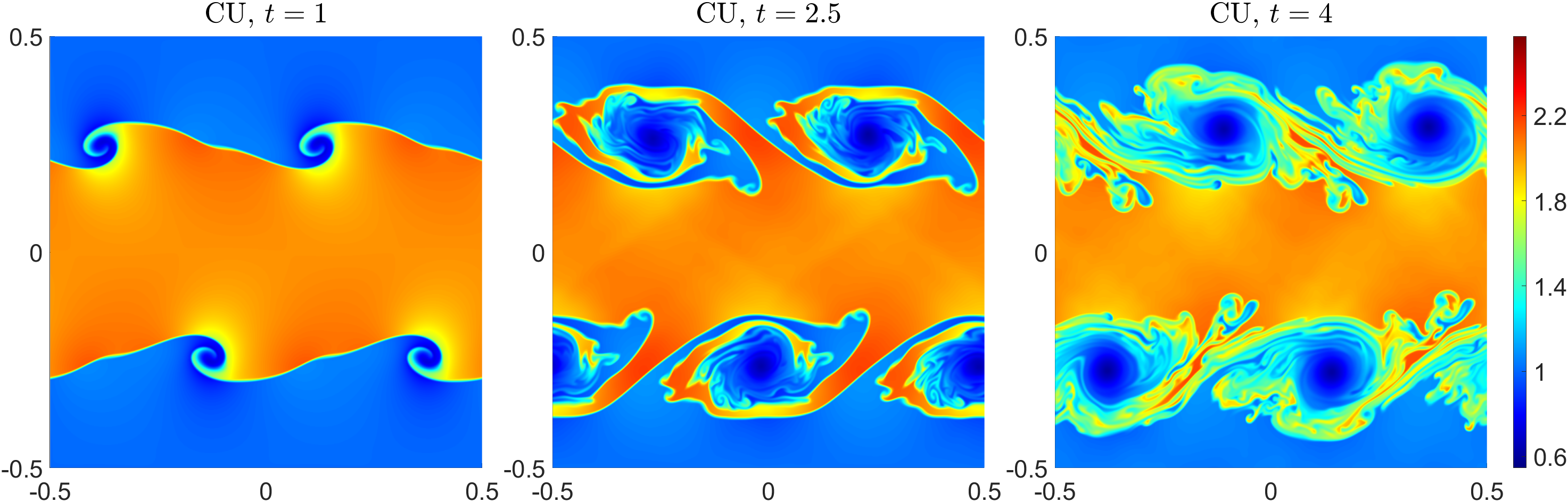}}
\vskip8pt
\centerline{\includegraphics[trim=0cm 0cm 0cm 0cm, clip, width=\linewidth]{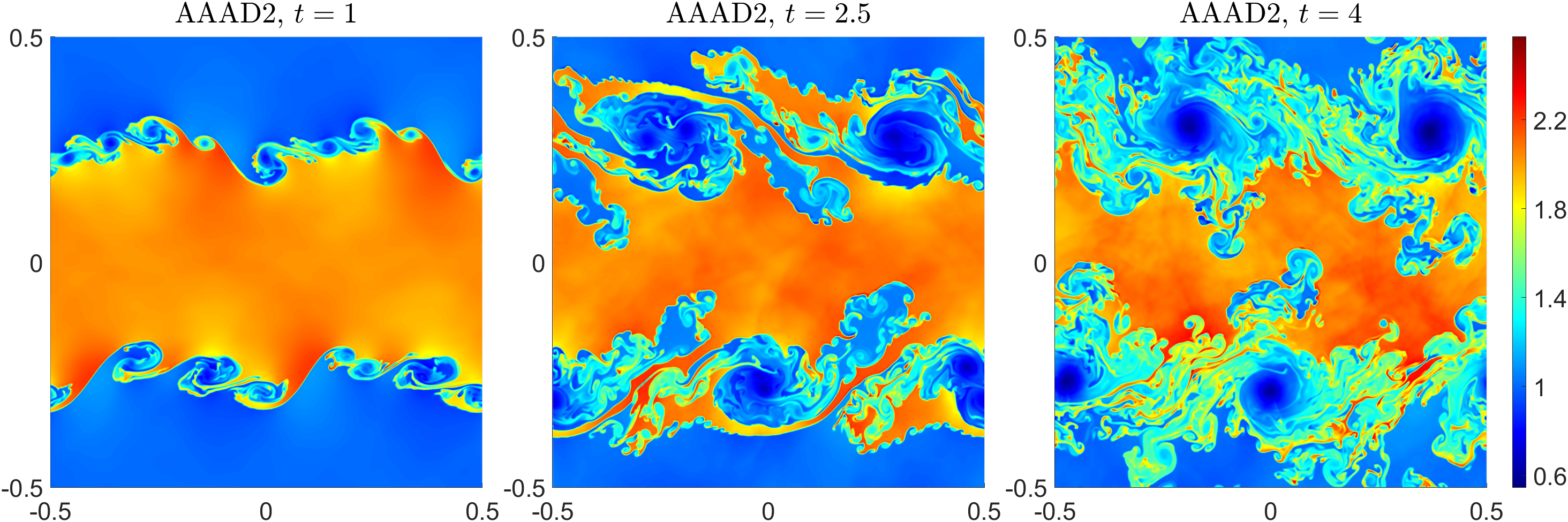}}
\caption{\sf Example 13: Density $\rho$ computed by the CU (top row) and AAAD2 (bottom row) schemes at $t=1$ (left column), $2.5$ (middle
column), and $4$ (right column).\label{fig12a}}
\end{figure}

\begin{figure}[ht!]
\centerline{\includegraphics[trim=0cm 0cm 0cm 0cm, clip, width=\linewidth]{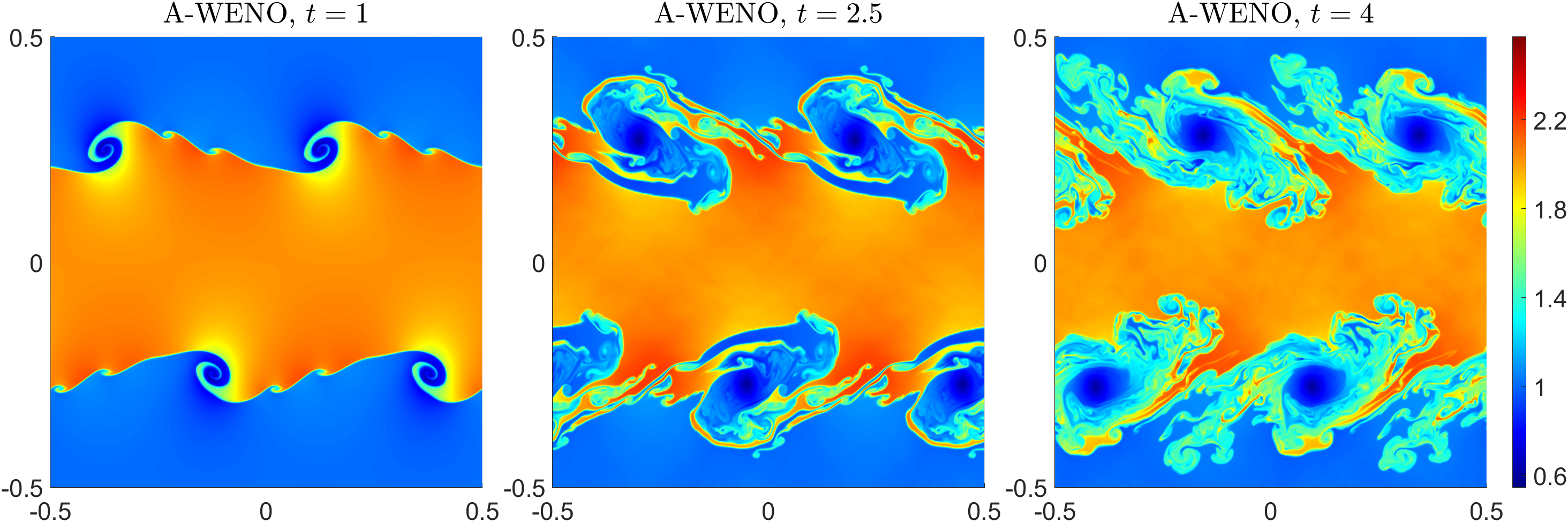}}
\vskip8pt
\centerline{\includegraphics[trim=0cm 0cm 0cm 0cm, clip, width=\linewidth]{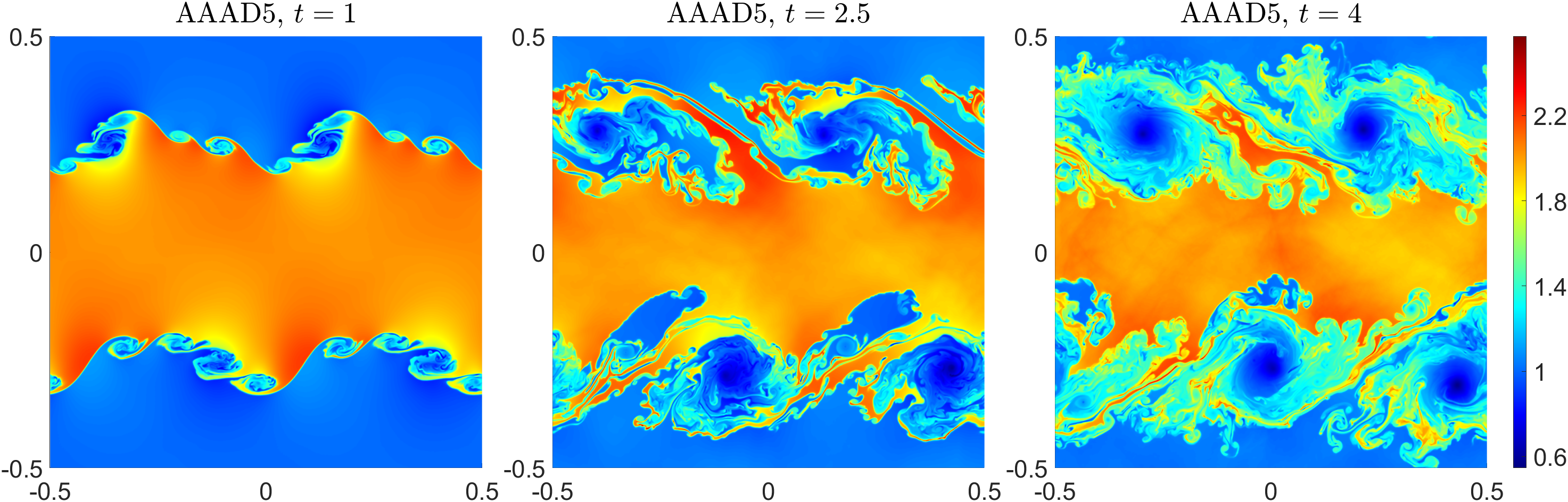}}
\caption{\sf Example 13: Density $\rho$ computed by the A-WENO (top row) and AAAD5 (bottom row) schemes at $t=1$ (left column), $2.5$
(middle column), and $4$ (right column).\label{fig12b}}
\end{figure}

We compute the numerical solution until the final time $t=4$ by the studied AAAD2 (with the adaptation constant $\texttt C=0.05$) and AAAD5
(with the adaptation constant $\texttt C=0.01$) schemes on a uniform mesh with $\dx=\dy=1/1024$ and plot the numerical results at times
$t=1$, $2.5$, and $4$ in Figures \ref{fig12a}--\ref{fig12b}. One can observe that already at $t=1$, the vortex streets produced by the AAAD
schemes are more pronounced. These structures grow over time, leading to increasingly complex mixing patterns at later times. Overall, the
AAAD schemes outperform their non-adaptive counterparts by capturing richer turbulent structures.

\subsubsection*{Example 14---Rayleigh-Taylor (RT) Instability}
In the last example taken from \cite{Shi03}, we investigate the RT instability, which occurs when a layer of heavier fluid is placed above a
layer of lighter fluid. The model is governed by the 2-D Euler equations with gravitational source terms acting in the positive
$y$-direction, namely,
\begin{equation*}
\begin{aligned}
&\rho_t+(\rho u)_x+(\rho v)_y=0,\\
&(\rho u)_t+(\rho u^2+p)_x+(\rho uv)_y=0,\\
&(\rho v)_t+(\rho uv)_x+(\rho v^2+p)_y=\rho,\\
&E_t+\left[u(E+p)\right]_x+\left[v(E+p)\right]_y=\rho v.
\end{aligned}
\end{equation*}
We consider the initial conditions
\begin{equation*}
(\rho(x,y,0),u(x,y,0),v(x,y,0),p(x,y,0))=\begin{cases}(2,0,-0.025\,c\cos(8\pi x),2y+1),&y<0.5,\\
(1,0,-0.025\,c\cos(8\pi x),y+1.5),&\mbox{otherwise},\end{cases}
\end{equation*}
where $c$ is the speed of sound. The computational domain is $[0,0.25]\times[0,1]$ with the solid wall boundary conditions imposed at $x=0$
and $x=0.25$, and the following Dirichlet boundary conditions imposed at the top and bottom boundaries:
\begin{equation*}
(\rho,u,v,p)\big|_{y=1}=(1,0,0,2.5),\quad(\rho,u,v,p)\big|_{y=0}=(2,0,0,1).
\end{equation*}

\begin{figure}[ht!]
\centerline{\includegraphics[trim=0cm 0cm 0cm 0cm, clip, width=\linewidth]{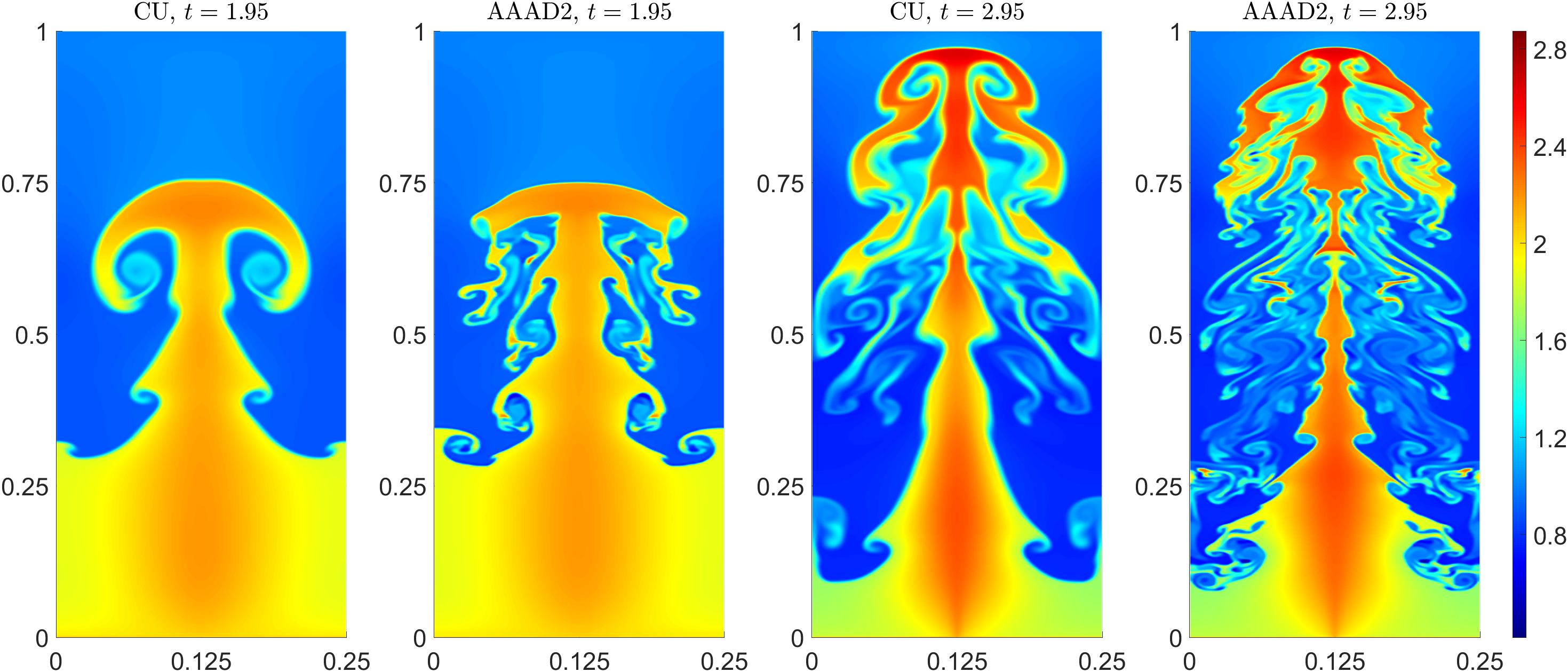}}
\vskip8pt
\centerline{\includegraphics[trim=0cm 0cm 0cm 0cm, clip, width=\linewidth]{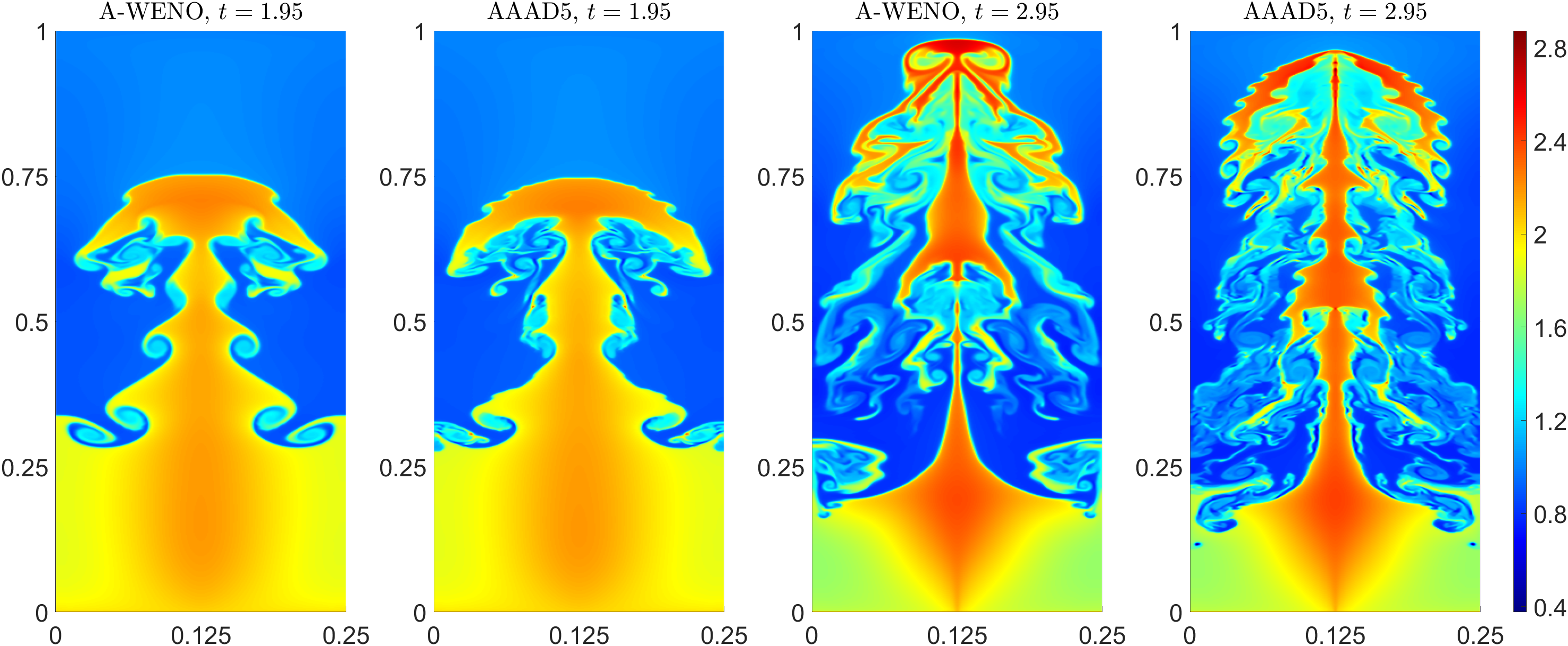}}
\caption{\sf Example 14: Density $\rho$ computed by the CU, AAAD2, A-WENO, and AAAD5 schemes at different times.\label{fig13}}
\end{figure}

We compute the numerical solution until the final time $t=2.95$ by the studied AAAD2 (with the adaptation constant $\texttt C=0.05$) and
AAAD5 (with the adaptation constant $\texttt C=0.02$) schemes on a uniform mesh with $\dx=\dy=1/1024$ and plot the results, obtained at
times $t=1.95$ and $2.95$, in Figure \ref{fig13}. One can observe pronounced differences between the AAAD solutions and those computed by
the CU and A-WENO schemes: the AAAD schemes capture substantially more complex small-scale structures, which again indicates reduced
numerical dissipation present in the AAAD schemes.

\section{Conclusion}\label{sec5}
In this paper, we have introduced a new adaptive artificial anti-diffusion (AAAD) methodology for substantially improving the resolution of
contact waves and complicated wave structures generated by multiple wave interactions, while preserving the stability and non-oscillatory
nature of the underlying numerical method for hyperbolic systems of conservation laws. The main idea is to enhance the resolution by
adaptively adding artificial numerical anti-diffusion (AD) in linearly degenerate characteristic fields only. In this way, the excessive
numerical dissipation present in the underlying numerical schemes can be reduced near contact waves, while genuinely nonlinear fields
associated with shocks and rarefactions are not affected by undesirable anti-diffusive corrections. As a result, the proposed AAAD methods
can produce highly accurate and essentially oscillation-free numerical solutions.

We have developed both second-order finite-volume and fifth-order finite-difference AAAD schemes. The AD coefficients are selected
adaptively with the help of smoothness indicators, which are used to distinguish between ``smooth'', ``rough'', and ``rough contact''
regions. This adaptive mechanism makes the AD sufficiently strong near contact discontinuities, where improved resolution is mostly needed,
while keeping it very small in ``smooth'' regions so that the formal order of accuracy of the underlying scheme is preserved.

Although the AAAD methods developed in this paper have been realized using central-upwind numerical fluxes, the proposed strategy is not
tied to any particular class of methods. Since the AD correction is added to a stable baseline solver and is constructed through the local
characteristic decomposition, the same idea can be directly combined with any robust and nonlinearly stable numerical methods for hyperbolic
systems of conservation laws. Therefore, the proposed AAAD approach can be viewed as a general resolution-enhancement procedure, which can
be used as an add-on mechanism for reducing excessive numerical dissipation in linearly degenerate fields.

The performance of the proposed methods has been demonstrated on a series of one- and two-dimensional numerical examples for the Euler
equations of gas dynamics. The obtained numerical results show that the AAAD schemes substantially sharpen contact discontinuities and
better resolve finer structures generated by complicated wave interactions, while maintaining robustness and suppressing spurious
oscillations. Although the present paper focuses on the Euler equations of gas dynamics, the proposed AAAD methodology can be naturally
extended to other hyperbolic systems, for which the characteristic structure is available and the linearly degenerate fields can be
explicitly identified.

\begin{DA}
\paragraph{Funding.} The work of S. Chu was funded by the Deutsche Forschungsgemeinschaft (DFG, German Research Foundation) - SPP 2410
Hyperbolic Balance Laws in Fluid Mechanics: Complexity, Scales, Randomness (CoScaRa) within the Project(s) HE5386/26-1 (Numerische
Verfahren für gekoppelte Mehrskalenprobleme,525842915) and (Zufällige kompressible Euler Gleichungen: Numerik und ihre Analysis, 525853336)
HE5386/27-1, and  the Deutsche Forschungsgemeinschaft (DFG, German Research Foundation) - SPP 2183: Eigenschaftsgeregelte Umformprozesse
with the Project(s) HE5386/19-2,19-3 Entwicklung eines flexiblen isothermen Reckschmiedeprozesses für die eigenschaftsgeregelte Herstellung
von Turbinenschaufeln aus Hochtemperaturwerkstoffen (424334423). The work of A. Kurganov was supported in part by NSFC grant W2431004.

\paragraph{Conflicts of interest.} On behalf of all authors, the corresponding author states that there is no conflict of interest.

\paragraph{Data and software availability.} The data that support the findings of this study as well as FORTRAN codes developed by the
authors and used to obtain all of the presented numerical results, are available from the corresponding author upon reasonable request.
\end{DA}

\appendix
\section{1-D Central-Upwind Numerical Fluxes}\label{appa}
The CU numerical fluxes introduced in \cite{Kurganov07} are
\begin{equation}
\bm{{\cal F}}^{\rm FV}_\jph=\frac{a^+_\jph\mF(\mU^-_\jph)-a^-_\jph\mF(\mU^+_\jph)}{a^+_\jph-a^-_\jph}+
\frac{a^+_\jph a^-_\jph}{a^+_\jph-a^-_\jph}\Big(\mU^+_\jph-\mU^-_\jph-\bm q_\jph\Big).
\label{2.2}
\end{equation}
Here, $\mU^\pm_\jph$ are one-sided point values of $\mU$ at the interface $x=x_\jph$ reconstructed from the cell averages $\{\,\xbar\mU_j\}$
using a nonlinear limiter applied to local characteristic variables.

In \eref{2.2}, $a^\pm_\jph$ are one-sided local propagation speeds, which can be estimated with the help of the smallest and largest
eigenvalues of the Jacobian $A(\mU)$. For example, one may take
$$
a^+_\jph=\max\Big\{u_\jph^-+c_\jph^-,u_\jph^++c_\jph^+,0\Big\},\quad
a^-_\jph=\min\Big\{u_\jph^--c_\jph^-,u_\jph^+-c_\jph^+,0\Big\},
$$
where
$$
u_\jph^\pm=\frac{(\rho u)_\jph^\pm}{\rho_\jph^\pm},\quad
p_\jph^\pm=(\gamma-1)\left[E_\jph^\pm-\hf\rho_\jph^\pm\big(u_\jph^\pm\big)^2\right],\quad
c_\jph^\pm=\sqrt{\frac{\gamma p_\jph^\pm}{\rho_\jph^\pm}}.
$$

Finally, $\bm q_\jph$ in \eref{2.2} is the ``built-in'' AD term (not to be confused with the AAAD term introduced in this paper), which was
derived in \cite{Kurganov07} and is given by
$$
\bm q_\jph={\rm minmod}\Big(\mU^+_\jph-\mU^*_\jph,\,\mU^*_\jph-\mU^-_\jph\Big)
$$
with
$$
\mU^*_\jph=\frac{a^+_\jph\mU^+_\jph-a^-_\jph\mU^-_\jph-\Big(\mF\big(\mU^+_\jph\big)-\mF\big(\mU^-_\jph\big)\Big)}{a^+_\jph-a^-_\jph}.
$$

\section{1-D LCD-Based MinMod2 Reconstruction}\label{appb}
Given the matrices $R_\jph$ and $\big(R_\jph\big)^{-1}$ in \eref{star3}, we first introduce the local characteristic variables $\bm\Gamma$
in the neighborhood of $x=x_\jph$:
$$
\bm\Gamma_\ell=\big(R_\jph\big)^{-1}\,\xbar\mU_{j+\ell},\quad\ell=-1,0,1,2.
$$
Equipped with the values $\bm\Gamma_{-1}$, $\bm\Gamma_0$, $\bm\Gamma_1$, and $\bm\Gamma_2$, we compute the slopes
\begin{equation}
{(\bm\Gamma_x)}_0={\rm minmod}\left(2\,\frac{\bm\Gamma_0-\bm\Gamma_{-1}}{\dx},\,\frac{\bm\Gamma_1-\bm\Gamma_{-1}}{2\dx},\,
2\,\frac{\bm\Gamma_1-\bm\Gamma_0}{\dx}\right),
\label{A1}
\end{equation}
and 
\begin{equation}
{(\bm\Gamma_x)}_1={\rm minmod}\left(2\,\frac{\bm\Gamma_1-\bm\Gamma_0}{\dx},\,\frac{\bm\Gamma_2-\bm\Gamma_0}{2\dx},\,
2\,\frac{\bm\Gamma_2-\bm\Gamma_1}{\dx}\right).
\label{A2}
\end{equation}
Here, the minmod function, defined as
\begin{equation}
{\rm minmod}(z_1,z_2,\ldots):=\begin{cases}
\min_j\{z_j\}&\mbox{if}~z_j>0\quad\forall\,j,\\
\max_j\{z_j\}&\mbox{if}~z_j<0\quad\forall\,j,\\
0            &\mbox{otherwise},
\end{cases}
\label{B3f}
\end{equation}
is applied in the component-wise manner.

Equipped with \eref{A1} and \eref{A2}, we evaluate
\begin{equation}
\bm\Gamma^-_\hf=\bm\Gamma_0+\frac{\dx}{2}({\bm\Gamma_x)}_0,\quad
\bm\Gamma^+_\hf=\bm\Gamma_1-\frac{\dx}{2}{(\bm\Gamma_x)}_1,
\label{A3}
\end{equation}
and then obtain the corresponding point values of $\mU$ by
\begin{equation}
\mU^\pm_\jph=R_\jph\bm\Gamma^\pm_\hf.
\label{B4}
\end{equation}

\section{1-D LCD-Based Fifth-Order WENO-Z Interpolant}\label{appc}
Given the point values $\psi_j$ of a certain function $\psi$ at uniform grid points $x=x_j$, the values $\psi^-_\jph$ are computed using a
weighted average of the three parabolic interpolants ${{\cal P}}_0(x)$, ${{\cal P}}_1(x)$ and ${{\cal P}}_2(x)$ obtained using the stencils
$[x_{j-2},x_{j-1},x_j]$, $[x_{j-1},x_j,x_{j+1}]$, and $[x_j,x_{j+1},x_{j+2}]$, respectively:
\begin{equation}
\psi^-_\jph=\sum_{k=0}^2\omega_k{\cal P}_k(x_\jph),
\label{C1}
\end{equation}
where
\begin{equation}
\begin{aligned}
&{\cal P}_0(x_\jph)=\frac{3}{8}\psi_{j-2}-\frac{5}{4}\psi_{j-1}+\frac{15}{8}\psi_j,\\
&{\cal P}_1(x_\jph)=-\frac{1}{8}\psi_{j-1}+\frac{3}{4}\psi_j+\frac{3}{8}\psi_{j+1},\\
&{\cal P}_2(x_\jph)=\frac{3}{8}\psi_j+\frac{3}{4}\psi_{j+1}-\frac{1}{8}\psi_{j+2}.
\end{aligned}
\label{C2}
\end{equation}
To ensure \eref{C1}--\eref{C2} are fifth-order accurate and nonoscillatory, we take the following weights $\omega_k$ in \eref{C1}:
\begin{equation}
\omega_k:=\frac{\alpha_k}{\alpha_0+\alpha_1+\alpha_2},\quad
\alpha_k=d_k\left[1+\left(\frac{\tau_5}{\beta_k+\varepsilon}\right)^p\right],\quad\tau_5=|\beta_2-\beta_0|,
\label{C3}
\end{equation}
where $d_0=\frac{1}{16}$, $d_1=\frac{5}{8}$, and $d_2=\frac{5}{16}$, and the SIs $\beta_k$ are given by 
\begin{equation}
\begin{aligned}
&\beta_0=\frac{13}{12}\big(\psi_{j-2}-2\psi_{j-1}+\psi_j\big)^2+\frac{1}{4}\big(\psi_{j-2}-4\psi_{j-1}+3\psi_j\big)^2,\\
&\beta_1=\frac{13}{12}\big(\psi_{j-1}-2\psi_j+\psi_{j+1}\big)^2+\frac{1}{4}\big(\psi_{j-1}-\psi_{j+1}\big)^2,\\
&\beta_2=\frac{13}{12}\big(\psi_j-2\psi_{j+1}+\psi_{j+2}\big)^2+\frac{1}{4}\big(3\psi_{j}-4\psi_{j+1}+\psi_{j+2}\big)^2.
\end{aligned}
\label{C4}
\end{equation}
In this paper, we have used $p=2$ and $\eps=10^{-12}$ in all of the numerical examples. The corresponding right-sided value, $\psi_\jph^+$,
can also be derived using a mirror-symmetric approach, and we omit it for the sake of brevity.  

To ensure a nonoscillatory nature of the reconstruction \eref{C1}--\eref{C4}, we apply it to the local characteristic variables, which we
introduce in the neighborhood of $x=x_\jph$:
$$
\bm\Gamma_\ell=\big(R_\jph\big)^{-1}\mU_{j+\ell},\quad\ell=-2,\dots,3.
$$
Equipped with the values $\bm\Gamma_{-2}$, $\bm\Gamma_{-1}$, $\bm\Gamma_0$, $\bm\Gamma_1$, $\bm\Gamma_2$, and $\bm\Gamma_3$, we then apply
the interpolation procedure \eref{C1}--\eref{C4} to each of the components $\Gamma^{(i)}$, $i=1,\dots,d$ of $\bm\Gamma$ to obtain
${\bm\Gamma}_\hf^-$ (the values of ${\bm\Gamma}_\hf^+$ are computed in the mirror-symmetric way). Finally, the corresponding one-sided point
values of $\mU$ are given by \eref{B4}.

\section{2-D Central-Upwind Numerical Fluxes}\label{appd}
The 2-D CU numerical fluxes from \cite{Kurganov07,CCKST2018} are
\begin{equation}
\hspace*{-0.1cm}\begin{aligned}
\bm{{\cal F}}^{\rm FV}_{\jph,k}&=\frac{a^+_{\jph,k}\mF\big(\mU^-_{\jph,k}\big)-a^-_{\jph,k}\mF\big(\mU^+_{\jph,k}\big)}
{a^+_{\jph,k}-a^-_{\jph,k}}+\frac{a^+_{\jph,k}a^-_{\jph,k}}{a^+_{\jph,k}-a^-_{\jph,k}}\left(\mU^+_{\jph,k}-\mU^-_{\jph,k}-\bm q^x_{\jph,k}
\right),\\[0.5ex]
\bm{{\cal G}}^{\rm FV}_{j,\kph}&=\frac{b^+_{j,\kph}\mG\big(\mU^-_{j,\kph}\big)-b^-_{j,\kph}\mG\big(\mU^+_{j,\kph}\big)}
{b^+_{j,\kph}-b^-_{j,\kph}}+\frac{b^+_{j,\kph}b^-_{j,\kph}}{b^+_{j,\kph}-b^-_{j,\kph}}\left(\mU^+_{j,\kph}-\mU^-_{j,\kph}-\bm q^y_{j,\kph}
\right).
\end{aligned}
\label{4.4a}
\end{equation}
Here, $\mU^\pm_{\jph,k}$ and $\mU^\pm_{j,\kph}$ are the one-sided interface values obtained at the points $(x,y)=(x_\jph,y_k)$ and
$(x,y)=(x_j,y_\kph)$, respectively, using a piecewise linear reconstruction applied to the local characteristic variables.

The one-sided local propagation speeds $a^\pm_{\jph,k}$ and $b^\pm_{j,\kph}$ can be estimated using the smallest and largest eigenvalues of
the Jacobians $A(\mU)$ and $B(\mU)$. We set
$$
\begin{aligned}
&a^+_{\jph,k}=\max\Big\{u^-_{\jph,k}+c^-_{\jph,k},u^+_{\jph,k}+c^+_{\jph,k},0\Big\},
&&a^-_{\jph,k}=\min\Big\{u^-_{\jph,k}-c^-_{\jph,k},u^+_{\jph,k}-c^+_{\jph,k},0\Big\},\\
&\,b^+_{j,\kph}=\max\Big\{v^-_{j,\kph}+c^-_{j,\kph},\,v^+_{j,\kph}+c^+_{j,\kph},\,0\Big\},
&&\,b^-_{j,\kph}=\min\Big\{v^-_{j,\kph}-c^-_{j,\kph},\,v^+_{j,\kph}-c^+_{j,\kph},\,0\Big\}.
\end{aligned}
$$
The ``built-in'' AD terms $\bm q^x_{\jph,k}$ and $\bm q^y_{j,\kph}$ are defined as in \cite{CCKST2018}:
$$
\begin{aligned}
\bm q^x_{\jph,k}&={\rm minmod}\Big(\mU^+_{\jph,k}-\mU^*_{\jph,k},\,\mU^*_{\jph,k}-\mU^-_{\jph,k}\Big),\\
\bm q^y_{j,\kph}&={\rm minmod}\Big(\mU^+_{j,\kph}-\mU^*_{j,\kph},\,\mU^*_{j,\kph}-\mU^-_{j,\kph}\Big),
\end{aligned}
$$
where
$$
\begin{aligned}
\mU^*_{\jph,k}&=\frac{a^+_{\jph,k}\mU^+_{\jph,k}-a^-_{\jph,k}\mU^-_{\jph,k}-
\Big(\mF\big(\mU^+_{\jph,k}\big)-\mF\big(\mU^-_{\jph,k}\big)\Big)}{a^+_{\jph,k}-a^-_{\jph,k}},\\
\mU^*_{j,\kph}&=\frac{b^+_{j,\kph}\mU^+_{j,\kph}-b^-_{j,\kph}\mU^-_{j,\kph}-
\Big(\mG\big(\mU^+_{j,\kph}\big)-\mG\big(\mU^-_{j,\kph}\big)\Big)}{b^+_{j,\kph}-b^-_{j,\kph}}.
\end{aligned}
$$

\section{2-D LCD-Based MinMod2 Reconstruction}\label{appe}
Here, we describe how the cell interface point values of $\bm U$ can be obtained using the 2-D MinMid2 reconstruction applied in the local
characteristic variables. We will only show how to compute the point values $\mU^\pm_{\jph,k}$, as $\mU^\pm_{j,\kph}$ can be calculated in a
similar way.

Given the matrices $R_{\jph,k}$ and $\big(R_{\jph,k}\big)^{-1}$, we first introduce the local characteristic variables
$\bm\Gamma$ in the neighborhood of $(x,y)=(x_\jph,y_k)$:
\begin{equation*}
\bm\Gamma_\ell=\big(R_{\jph,k}\big)^{-1}\,\xbar\mU_{j+\ell,k},\quad\ell=-1,0,1,2.
\end{equation*}
Equipped with the values $\bm\Gamma_{-1}$, $\bm\Gamma_0$, $\bm\Gamma_1$, and $\bm\Gamma_2$, we compute $\bm\Gamma^\pm_\hf$ precisely as in
\eref{A1}--\eref{A3}, and then 
\begin{equation*}
\mU^\pm_{\jph,k}=R_{\jph,k}\bm\Gamma^\pm_\hf.
\end{equation*}

\bibliographystyle{siamnodash}
\bibliography{ref}
\end{document}